\documentclass[journal]{IEEEtran}

\usepackage{graphicx}
\usepackage[colorlinks=true,
            linkcolor=blue,
            citecolor=blue,
            urlcolor=blue]{hyperref}
\usepackage{psfrag}
\usepackage{hyphenat}
\usepackage{relsize} 
\usepackage{amsmath}
\usepackage{amsfonts}
\usepackage[dvipsnames]{xcolor}
\usepackage{bm}
\usepackage{paralist}
\usepackage[caption=false]{subfig}
\usepackage{comment}
\usepackage{tabularx}
\usepackage{algorithm}
\usepackage{algpseudocode}
\algrenewcommand\algorithmiccomment[1]{\hfill\(\triangleright\) #1}
\usepackage{siunitx}
\usepackage{booktabs}
\usepackage{multirow}
\usepackage{rotating}
\usepackage{pdflscape}
\usepackage{array}
\usepackage{tabularx}
\usepackage[table]{xcolor}    
\usepackage{makecell}
\usepackage{cite}
\newcolumntype{Y}{>{\raggedright\arraybackslash}X}  

\newcommand{\mvar}[2]{#1_{\mathrm{#2}}}

\begin{document}
\setlength{\abovedisplayskip}{3pt}
\setlength{\belowdisplayskip}{3pt}

\setlength{\textfloatsep}{4pt plus 1.0pt minus 2.0pt}

\setlength{\medmuskip}{0mu}
\setlength{\thickmuskip}{0mu}
\setlength{\thinmuskip}{0mu}
\title{PriEco-DRL: Joint Optimization of  Electric-Bus Eco-Driving and Transit-Priority Adaptive Signals via Deep Reinforcement Learning}

\author{Dingshan Sun, Ang Li, Chaopeng Tan, Zicheng Su, Wanjing Ma, and Marco Rinaldi%
\thanks{*Corresponding author: Chaopeng Tan (email: chaopeng.tan@tu-dresden.de).}

\thanks{Dingshan Sun, Zicheng Su and Wanjing Ma are with the Key Laboratory of Road and Traffic Engineering of the Ministry of Education, College of Transportation, Tongji University, Shanghai 201804, China.}%
\thanks{Ang Li is with the Department of Aeronautical and Aviation Engineering, The Hong Kong Polytechnic University, Kowloon, Hong Kong, China.}%
\thanks{Chaopeng Tan is with the Chair of Traffic Process Automation, Technische Universit\"at Dresden, 01069 Dresden, Germany.}%
\thanks{Marco Rinaldi is with the Department of Transport and Planning, Delft University of Technology, 2628 CN Delft, The Netherlands.}%
\thanks{The authors would like to acknowledge partial financial support of the National Natural Science Foundation of China (No. 52325210) and the EU Horizon Europe Research and Innovation Programme
(Grant Agreement No. 101103808 ACUMEN).}%
}
\maketitle
\begin{abstract}
Urban transit electrification requires balancing energy efficiency, schedule reliability, and ride comfort for electric buses (EBs), particularly when interacting with transit-priority adaptive signals in congested networks. This paper proposes PriEco-DRL, a joint optimization framework that integrates EB eco-driving with transit-priority adaptive signal control using deep reinforcement learning (DRL). The signal layer employs a priority-weighted max-pressure (Priority-MP) controller to allocate green time based on occupancy-aware pressures, while the vehicle layer adapts longitudinal control based on local uncertain, dynamically evolving signal cues.
A structured reward combines guidance and event-based reinforcement to align EB arrivals with green opportunities, while considering energy, time, comfort, and safety. The framework uses centralized training and decentralized execution (CTDE) with parameter sharing, allowing a single DRL agent to learn from multiple buses and routes using local observations.
Experiments on a real-world corridor show that PriEco-DRL reduces EB energy consumption while maintaining network efficiency and transit priority compared to fixed-time, actuated, and rule-based signal–vehicle coordination baselines. Energy- and trajectory-based analyses reveal that the improvements stem from fewer unscheduled stop–start events and smoother speed regulation under adaptive signals. The results highlight a tunable energy–time trade-off, allowing flexible operational choices through reward weighting.

\end{abstract}

\begin{IEEEkeywords}
Eco-driving; Transit priority signal; Max-pressure; Deep reinforcement learning;
Electric bus.
\end{IEEEkeywords}
\section{Introduction}
\IEEEPARstart{W}{ith} ongoing urbanization, cities are becoming increasingly crowded, placing substantial pressure on urban transportation systems. Meanwhile, vehicle emissions continue to degrade the urban living environment. Promoting public transportation and reducing private car use is therefore widely recognized as an effective way to alleviate both congestion and environmental impacts \cite{banister2008sustainable,miller2016public}. Within modern public transport systems, electric buses (EBs) have become an important component due to their environmental benefits and their potential to support sustainable urban mobility \cite{manzolli2022review}.
However, EB operation also faces significant challenges. A key issue is how to jointly balance operational efficiency, operating cost, and passenger experience under complex and highly variable driving conditions, including fluctuating traffic density, temperature, passenger occupancy, and terrain \cite{papa2022electric,alarrouqi2024electric}. In addition, range anxiety and battery degradation remain important concerns, especially under repeated charging--discharging cycles and frequent stop-and-go driving \cite{rainieri2023psychological,zeng2022role}. Developing advanced EB operation strategies is therefore of great practical importance.

Existing studies on EB operation can be broadly divided into macro-level and micro-level approaches \cite{zhang2023vehicle}. Macro-level studies mainly focus on long-term planning, such as route optimization and charging scheduling \cite{zhang2024routing,zhou2024charging}. Although effective for strategic planning, such methods rely on historical data and cannot adapt to dynamically changing traffic conditions or unexpected events. In contrast, micro-level approaches provide real-time speed guidance, i.e., eco-driving, to reduce stop-and-go behavior and improve energy efficiency \cite{ma2013rule,zhang2024eco,pan2025eco}. By leveraging signal phase and timing (SPaT) information and vehicle-to-infrastructure (V2I) communication, these methods can align bus arrivals with feasible green windows and generate speed profiles that account for safety, comfort, and energy use.

However, most existing eco-driving studies assume fixed or deterministic traffic signals \cite{li2023overcome,yang2024eco,wegener2021automated}. In practice, adaptive signal control has become increasingly important in intelligent transportation systems, as it adjusts phase selection and green time in real time according to traffic demand and system performance. Compared with fixed-time control, adaptive signals can better reduce delay and spillback and support transit-priority operations \cite{lin2015transit}. 
Yet this also makes SPaT stochastic and time-varying, creating a challenging environment for conventional eco-driving strategies that rely on predictable green windows. This is particularly problematic for bus operations, which are inherently multi-stage due to fixed routes, designated stops, dwell, and post-stop departures. Methods developed for simplified or single-stage scenarios therefore transfer poorly to realistic bus settings. For EBs, the challenge is further compounded by the need to explicitly balance energy consumption with travel efficiency and service reliability. The key difficulty is therefore eco-driving under uncertain, dynamically evolving SPaT rather than under deterministic future green windows. To the best of our knowledge, few studies have explicitly addressed EB eco-driving under adaptive or transit-priority signal control \cite{yang2024entire}.

To address these gaps, we propose PriEco-DRL, a coupled framework that integrates priority-weighted max-pressure (Priority-MP) adaptive signal control with a deep reinforcement learning (DRL) eco-driving policy for electric buses. Priority-MP reallocates green time using occupancy-aware pressures to improve network efficiency while preserving transit priority. On the vehicle side, a DRL agent learns longitudinal speed guidance under dynamically evolving, pressure-driven SPaT. The policy operates over complete EB routes and explicitly accounts for multi-stage driving contexts, including cruise, stop approach, dwell, and intersection crossing. By modeling the signal--vehicle coupling explicitly, PriEco-DRL aims to reduce energy consumption while maintaining service efficiency and operational smoothness under uncertain traffic conditions and concurrent multi-route bus operations. The main contributions of this paper are as follows:
\begin{itemize}
    \item \textit{Coupled priority--eco-driving framework.} We propose PriEco-DRL, which explicitly couples a Priority-MP controller with a DRL eco-driving policy for electric buses. The policy uses local signal and pressure information to coordinate speed decisions with uncertain adaptive transit-priority signal control, extending EB eco-driving beyond fixed-time settings.

    \item \textit{Generalizable training and deployment regime.} We adopt centralized training, decentralized execution (CTDE) with parameter sharing and a network-agnostic state space. This improves data efficiency during training and supports decentralized deployment across multiple routes, layouts, and simultaneous EB operations.

    \item \textit{Real-world network evaluation.} We validate the proposed framework on a real Amsterdam arterial corridor with three signalized intersections, eight bus lines, and 18 stops. The policy is trained on a subset of routes and tested on unseen routes, demonstrating transferability under realistic mixed-lane and dedicated-lane operations.
\end{itemize}

The remainder of this paper is organized as follows. Section~\ref{sec:related-work} reviews the related work on adaptive signal control and eco-driving strategies. Section~\ref{sec:problem-statement} formulates the problem addressed in this study. Section~\ref{sec:MP} presents the Priority-MP signal control method. Section~\ref{sec:method} introduces the proposed PriEco-DRL framework. Section~\ref{sec:case-study} presents the case study and experimental results. Finally, Section~\ref{sec:conclusion} concludes the paper and outlines future research directions.

\section{Related work} \label{sec:related-work}
\subsection{Adaptive signal control}

Adaptive signal control has long been studied as an effective strategy for reducing delay, queues, and stops under time-varying demand. Early responsive methods, such as SCOOT \cite{hunt1981scoot} and OPAC \cite{gartner1983opac}, adjust signal timings online based on detector data and short-horizon traffic estimates, but their largely reactive nature limits their ability to address complex network interactions and rapidly changing traffic conditions. More advanced approaches, including model predictive control (MPC) and DRL, have therefore been increasingly explored. MPC enables prediction-based and constraint-aware optimization of traffic signals \cite{garcia1989model,ye2019survey,GUO2026103411}, while DRL learns adaptive control policies directly from interaction data without requiring an explicit predictive model \cite{haydari2020deep}. However, MPC often faces challenges in model fidelity and online computational tractability for large-scale urban networks, whereas DRL-based policies are typically less interpretable and less transparent for coordination with vehicle-side eco-driving decisions \cite{dulac2021challenges}. These limitations motivate the need for a signal control framework that is adaptive, computationally efficient, and sufficiently interpretable for integration with vehicle-side control.

Among real-time decentralized approaches, Max-Pressure (MP) control has emerged as a particularly promising strategy due to its simplicity, low computational burden, and provable stability properties under standard stochastic traffic assumptions \cite{varaiya2013max}. At each decision step, MP selects the phase associated with the highest upstream--downstream queue imbalance, thereby alleviating local bottlenecks while implicitly promoting network-wide coordination. MP requires only adjacent-link traffic states and does not depend on an explicit predictive model, making it well suited for real-time applications \cite{li2019position, liu2022novel}. In addition, its pressure formulation can be naturally extended to account for transit priority by weighting movements according to bus occupancy or delay \cite{xu2022integrating,tan2026cv,tan2025transit}. Compared with MPC- and DRL-based signal control, MP offers a more interpretable and computationally efficient infrastructure-side controller, which makes it particularly suitable for coupling with eco-driving strategies. Therefore, we adopt MP as the signal control backbone in this study. Specifically, we build on our recent Priority-MP design, which incorporates transit priority into pressure-based control while preserving the decentralized mechanism. Implementation details are provided in Section~\ref{sec:MP}.

\subsection{Eco-driving strategies}

With the emergence of autonomous driving and electric buses (EBs), eco-driving has become an important technology for improving energy efficiency while accounting for regenerative braking, auxiliary loads, and service reliability \cite{huang2018eco,mensing2014eco}. In general, eco-driving seeks longitudinal speed control strategies that balance energy consumption, traffic efficiency, and passenger comfort under vehicle, traffic, and signal constraints. Existing methods can be broadly categorized into three groups: rule-based heuristics, optimization-based control, and learning-based policies.

Rule-based approaches prescribe speed profiles through hand-crafted logic using local roadway and SPaT information. A representative family is GLOSA, which computes a target speed so that vehicles arrive within a feasible green window, thereby reducing stops and traction energy \cite{barth2011dynamic,eckhoff2013potentials}. For transit operations, GLODTA-type methods further coordinate dwell adjustment with downstream signal timing to improve green-wave alignment \cite{giorgione2017experimental,laskaris2020sensitivity}. These methods are attractive because of their simplicity, interpretability, and lightweight implementation. However, they usually rely on fixed thresholds and deterministic trigger logic, making them sensitive to calibration and less robust under stochastic traffic conditions and complex multi-objective trade-offs.

Optimization-based methods formulate eco-driving as a constrained optimal control problem, often solved by dynamic programming or MPC \cite{guo2016optimal,li2021online}. With SPaT/V2I information, these methods can explicitly encode signal timing, speed, jerk, and energy constraints, and are therefore well suited to energy-efficient approach and departure control \cite{yang2020eco,zhang2011predictive}. For EBs, optimization models can additionally incorporate regenerative braking, auxiliary loads, passenger-load-dependent mass, grade, and service-related objectives such as schedule adherence \cite{chow2017multi,jia2023novel}. Their main advantage lies in explicit constraint handling and physical interpretability. However, their effectiveness depends on model fidelity and reliable forecasts, while online computation becomes increasingly challenging in uncertain and highly dynamic traffic environments.

Learning-based approaches, especially DRL, provide an alternative by learning control policies directly from interaction data without relying on explicit predictive models \cite{li2022deep,wu2023deep}. DRL is particularly attractive for nonlinear, high-dimensional settings and can naturally optimize long-horizon objectives while incorporating safety or comfort constraints through reward shaping and action restrictions. Although DRL has been widely explored for vehicle eco-driving, studies on EB eco-driving remain limited. Existing works either focus on only part of the bus operation process or assume fixed/predictable signal timing \cite{huang2025towards,yang2024entire}, which restricts their applicability to real-world adaptive signal environments.

A smaller but growing body of literature considers eco-driving jointly with dynamic signal control. One line of work predicts vehicle arrivals and then adapts green windows to facilitate non-stop passage \cite{ma2013rule,sun2022eco}. Another line formulates a bi-level optimization problem, where the upper level optimizes signal timing and the lower level computes energy-efficient vehicle trajectories \cite{jung2016bi,li2018eco,ding2024intersection}. While effective in controlled settings, these methods remain strongly model-dependent and can become computationally demanding in complex corridors with mixed traffic and uncertain SPaT evolution. More importantly, existing formulations mainly target standard vehicles and do not explicitly account for EB-specific multi-stage operations such as stop approach, dwell, and post-stop departure.

In summary, the integration of EB eco-driving with adaptive and transit-priority signal control remains insufficiently explored. This paper addresses this gap by proposing PriEco-DRL, a coupled framework that combines a Priority-MP signal controller with a DRL-based eco-driving policy. The Priority-MP controller handles adaptive, occupancy-aware transit-priority signal control, while pressure-related signal information is provided to the DRL agent to support informed speed decisions under complex traffic conditions. In parallel, the reward design encodes multi-stage bus operations and multi-objective trade-offs, resulting in a practical eco-driving strategy for complete EB routes under adaptive signals.


\section{Problem statement}\label{sec:problem-statement}

\begin{figure}[!t]
    \centering
    \includegraphics[width=0.95\linewidth,trim={0cm, 3.5cm, 0cm, 3.5cm}, clip]{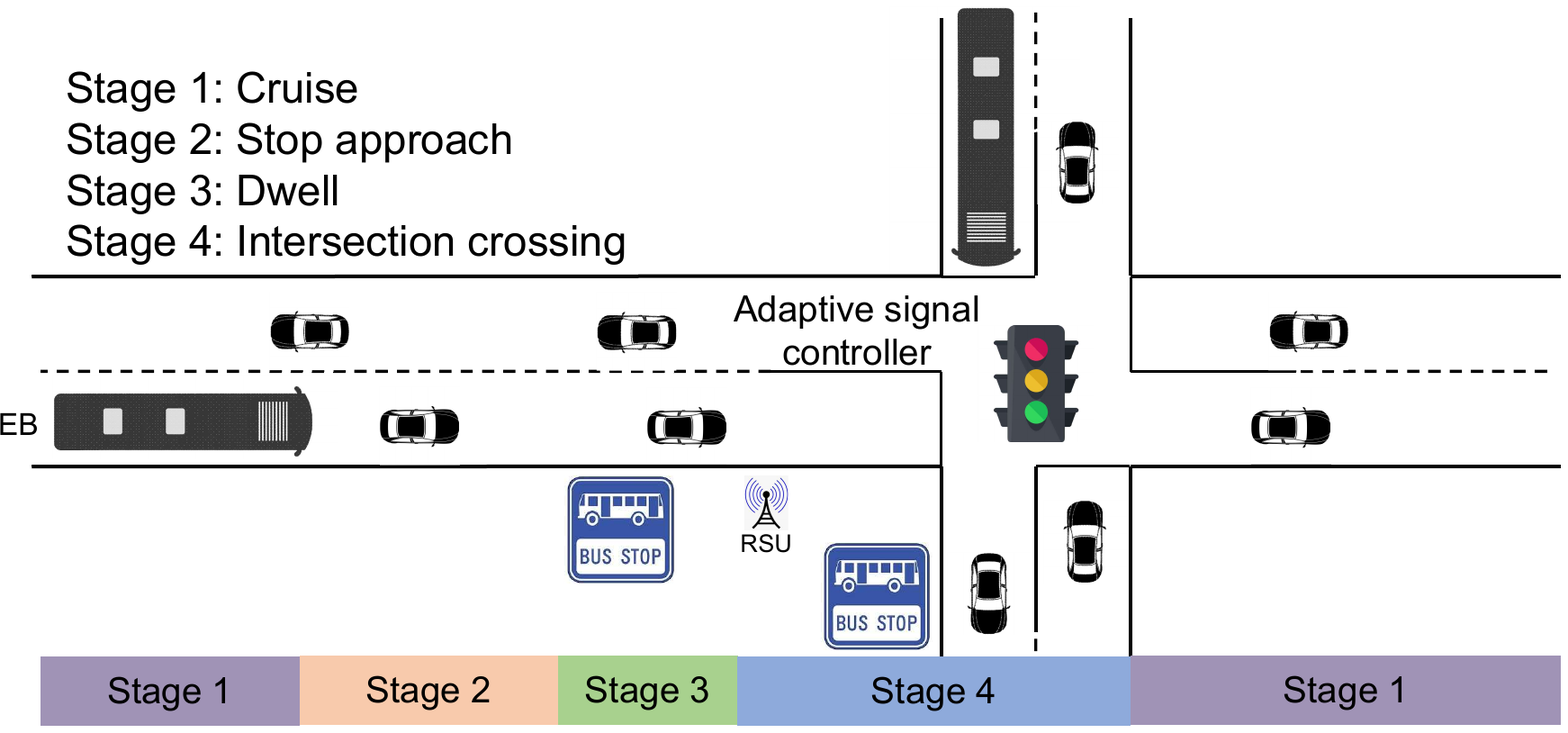}
    \caption{Stage division of EB entire route}
    \label{fig:Stages}
\end{figure}

In this paper, we consider EB eco-driving over an entire route. According to the driving context and operational objective, the route is divided into four stages: \emph{cruise}, \emph{stop approach}, \emph{dwell}, and \emph{intersection crossing}, as illustrated in Fig.~\ref{fig:Stages}. In the cruise stage, the EB travels between stops and intersections while maintaining safe car-following and energy-efficient motion. When approaching a stop, the EB enters the stop-approach stage, where it must decelerate safely and dock accurately. During the dwell stage, the EB remains at the stop for passenger boarding and alighting. In the intersection-crossing stage, the EB seeks to pass the intersection efficiently while avoiding unnecessary stopping and excessive energy use.

\begin{figure}[]
    \centering
    \includegraphics[width=1.0\linewidth, trim={0.2cm, 4.5cm, 0.1cm, 5.0cm}, clip]{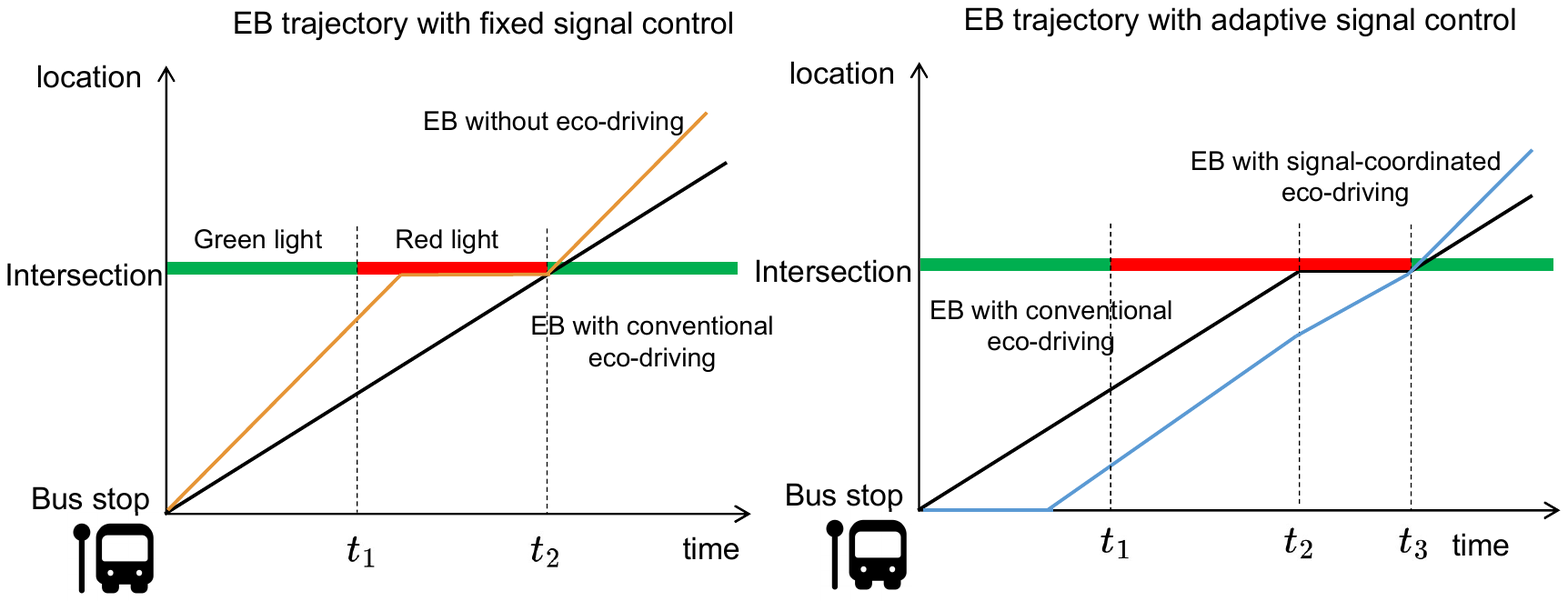}
    \vspace{-0.5cm}
    \caption{Comparisons of different driving strategies under different signal control methods}
    \label{fig:eco-driving}
\end{figure}

The main challenge lies in the intersection-crossing stage under adaptive signal control. As illustrated in Fig.~\ref{fig:eco-driving}, conventional eco-driving strategies typically compute a target speed by aligning the bus arrival time with an expected green window. This approach works when signal timing is fixed or predictable, but it becomes unreliable when the signal is adaptively adjusted according to real-time traffic conditions. In such cases, the green or red duration may be extended or terminated depending on competing traffic demands, so a bus following a conventional green-window strategy may still encounter an unexpected stop. This motivates a signal-coordinated eco-driving strategy that can exploit the information driving adaptive signal control and adjust bus speed accordingly.


The goal of this paper is therefore to develop a comprehensive EB eco-driving strategy that operates over the entire route and adapts to dynamic signal evolution. The overall objective is to jointly reduce trip energy consumption and travel time:
\begin{equation}
    \min(\omega_e\cdot\mvar{E}{trip}+\omega_t\cdot\mvar{T}{trip}), \label{eq:objecitve}
\end{equation}
where $\mvar{E}{trip}$ and $\mvar{T}{trip}$ denote the total energy consumption and travel time of a trip, and $\omega_e$ and $\omega_t$ are the corresponding weights. We assume that the required vehicle--infrastructure information exchange is available through connected infrastructure (e.g., RSUs/i-VICS). The main notations used in this paper are summarized in Table~\ref{tab:notations}.

\begin{table}[!t]
  \centering
  \caption{Notations of important variables}
  \label{tab:notations}
  \begin{tabularx}{\linewidth}{@{}lY@{}}
    \toprule
    \textbf{Variables} & \textbf{Explanation} \\
    \midrule

    \rowcolor{gray!10}\multicolumn{2}{@{}l}{\textbf{(A) Signal control variables}}\\
    $\bm{N}$      & The set of intersection nodes $n$ in the network \\
    $\bm{M}_n$      & The set of movements for a node $n$ \\
    $s_{(i,o)}$  & Signal state of movement $(i,o)\in \bm{M}_n$ where $i,o$ are link indices \\
    $\bm{s}_n$  & Signal vector of node $n$ composed of $s_{(i,o)}$\\
    $\bm{s}$ & The overall signal vector of the network composed of $\bm{s}_n$ \\
    $\mathcal{S}$ & The feasible space of all signal states \\
    $c_{(i,o)}$ & Saturation flow rate of movement $(i,o)$ \\
    $r_{(i,o)}$    & Turning ratio of movement $(i,o)$: the fraction of vehicles on incoming link $i$ to outgoing link $o$\\
    $u_j$  & The state of $j$-th vehicle in the movement \\
    $p_j$  & The number of passengers on the $j$-th vehicle in the movement \\
    $\alpha_j$  & The binary indicator that whether the $j$-th vehicle is counted for pressure calculation \\
    $x_j$  & The position of the $j$-th vehicle in the movement \\
    $\bm{J}_{(i,o)}$  & The set of all vehicles in the movement $(i,o)$ \\
    $P_{(i,o)}$ & The pressure of the movement $(i,o)$ \\
    \addlinespace[2pt]
    \midrule

    \rowcolor{gray!10}\multicolumn{2}{@{}l}{\textbf{(B) Eco-driving variables}}\\
    $v_t$  & The speed of ego EB at time step $t$  \\
    $a_t$ & The acceleration of ego EB at time step $t$\\
    $\sigma_t$ & The binary indicator of whether the EB is in the dwell stage at time step $t$\\
    $\alpha_t$ & The binary indicator of whether the EB has departed from the nearest upstream stop of the downstream intersection at time step $t$\\
    $\delta_t$ & Traffic signal status in the lane of EB (0, 1, 2 corresponds to green, yellow, and red respectively)\\
    $\mvar{\tau}{rem}$ & Remaining time for currently active signal phase\\
    $\bm{P}_t$ & Traffic signal pressure information vector: $\bm{P}_t=[P_{\phi_c}, P_{\phi_m}, \Delta P]$ where $P_{\phi_c}$ and $P_{\phi_m}$ are the aggregated pressure of specific phases $\phi_c$ and $\phi_m$, and $\Delta P$ is the difference between these two \\
    $\phi_c$ & The phase that includes the current movement of EB \\
    $\phi_m$ & The phase with the largest pressure among all other phases excluding $\phi_c$\\
    $a^{\mathrm{act}}_t$ & The action of DRL agent at time step $t$: the acceleration implemented on EB \\
    $d_l,d_{NS},d_{NI}$ & The distances from the EB to the leading vehicle, the next stop, and the stop line, respectively \\
    $d_L$ & The distance threshold beyond which the bus can be classified as initiating the stop maneuver at the next berth \\
    $D_{\mathrm{prep}}$ & The distance to identify whether the EB starts preparing to cross the intersection \\    
    $\underline{d}_{\mathrm{isec}},\bar{d}_{\mathrm{isec}}$ & The lower and upper bounds of $D_{\mathrm{prep}}$, respectively \\
    $T_c$ & Signal control step interval \\
    $\mvar{T}{safe},\mvar{T}{hard}$ & The safe and hard time thresholds of time to collision, respectively \\
    $\mvar{v}{entry},\mvar{V}{safe}$ & The bus stop entry speed of EB, and the maximum acceptable entry speed in the stop area  \\
    $\mathbb{E}[\Delta G]$ & The expected green extension of the current green phase\\
    \bottomrule
  \end{tabularx}
\end{table}

\section{Max-Pressure Control with Transit Priority} \label{sec:MP}

Max-Pressure (MP) control is an adaptive, model-free, and decentralized signal control strategy that selects, at each decision instant, the feasible phase plan with the largest upstream--downstream queue imbalance \cite{varaiya2013max}. In this paper, we adopt a stop-aware, occupancy-weighted Priority-MP controller developed in our previous work, and briefly summarize only the elements needed for the subsequent DRL formulation; full details can be found in \cite{tan2025transit}.

\begin{figure} 
\centering \includegraphics[width=0.48\textwidth]{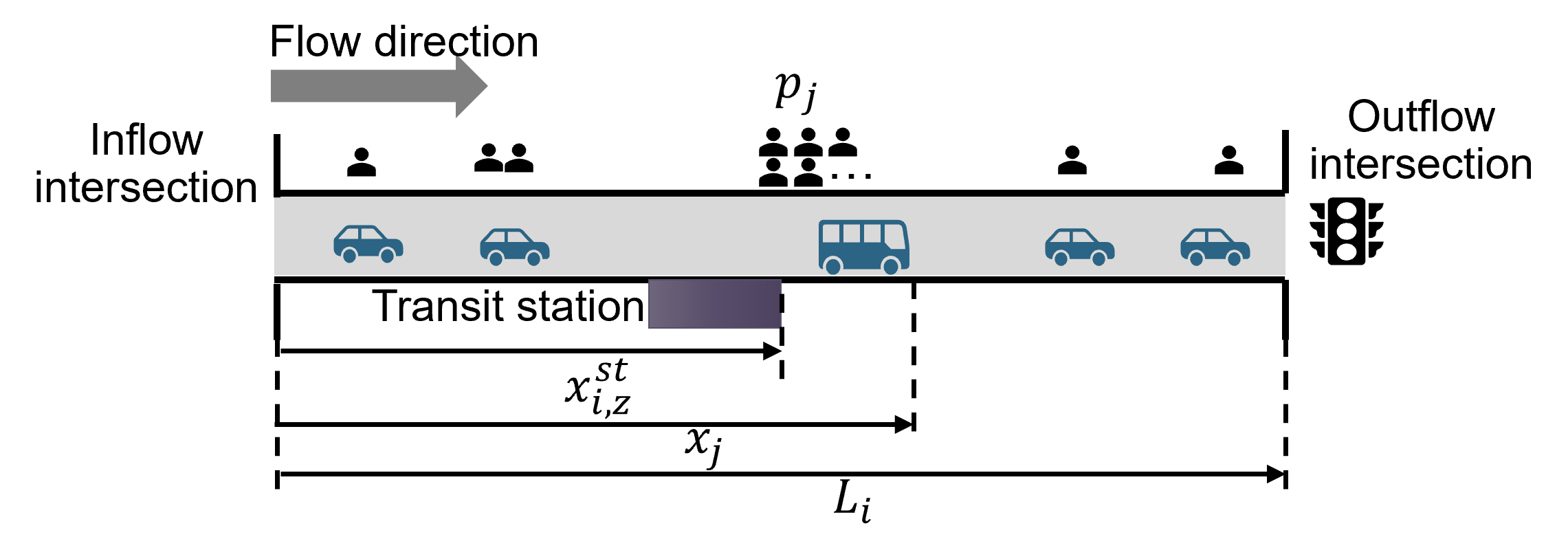} 
\caption{Parameters for Priority-MP considering transit stations} \label{fig: considerstop} 
\end{figure}

The key idea is to weight vehicle contributions by occupancy so that movements serving high-occupancy transit vehicles receive higher priority, as illustrated in Fig~\ref{fig: considerstop}. In addition, to avoid wasting green time when a bus is still dwelling at the nearest upstream stop, a transit vehicle contributes to pressure only after it has departed from that stop. The resulting movement pressure is defined as
\begin{align}
P_{(i,o)} ={}& s_{(i,o)} c_{(i,o)} \Biggl(
\sum_{j \in \bm{J}^{C}_{(i,o)}} p_j u_j
+ \sum_{j \in \bm{J}^{Tr}_{(i,o)}} \alpha_j p_j u_j \notag\\
&
- \sum_{(o,k)\in \bm{M}_{n'}} r_{(o,k)}
\Bigl(
\sum_{j \in \bm{J}^{C}_{(o,k)}} p_j u_j
+ \sum_{j \in \bm{J}^{Tr}_{(o,k)}} \alpha_j p_j u_j
\Bigr)
\Biggr), \label{eq:pressure-3}
\end{align}
and the signal plan is selected as
\begin{align}
\bm{s}^* ={}& \arg\max_{\bm{s}\in\mathcal{S}}
\sum_{n\in\bm{N}}
\Biggl(
\sum_{(i,o)\in\bm{M}_n} P_{(i,o)}
\Biggr). \label{eq:priority_MP_with_ts}
\end{align}
Here, $\bm{J}^{C}_{(i,o)}$ and $\bm{J}^{Tr}_{(i,o)}$ denote the sets of private cars and transit vehicles in movement $(i,o)$, respectively; $p_j$ is vehicle occupancy; $c_{(i,o)}$ is the saturation flow rate; and $r_{(o,k)}$ is the turning ratio toward downstream movement $(o,k)$.

The binary variable $\alpha_j$ determines whether a transit vehicle is counted in the pressure calculation:
\begin{align}
\alpha_j =
\begin{cases}
0, & \text{if } x_j < x_{i,1}^{st} \text{ or } T^{dwell}_j > 0,\\
1, & \text{otherwise,}
\end{cases}
\label{eq: alpha_j}
\end{align}
where $x_j$ is the vehicle position, $x_{i,1}^{st}$ is the position of the stop nearest to the downstream intersection on link $i$, and $T^{dwell}_j$ is the remaining planned dwell time.

The vehicle state used in the pressure calculation is the normalized link travel time,
\begin{equation}
u_j(t)=\frac{t-t^0_j}{\bar{\tau}_{(i,o)}},
\label{eq:state}
\end{equation}
where $t^0_j$ is the link-entry time of vehicle $j$ and $\bar{\tau}_{(i,o)}$ is the free-flow travel time of movement $(i,o)$. Eqs.~\eqref{eq:pressure-3}--\eqref{eq:state} define the Priority-MP controller used in this study. Note that we assume that vehicle link travel times can be obtained in real time through commonly available sensing sources such as dedicated short-range communications, automated vehicle identification, or LiDAR detectors. For more information on how detected vehicle penetration rates affect MP control, please refer to the study by \cite{tan2025transit, tan2026cv}.

\section{Deep reinforcement learning model for EB eco-driving} \label{sec:method}

This section presents the proposed PriEco-DRL framework for EB eco-driving. A DRL agent interacts with both the traffic signal controller and the traffic environment to learn a longitudinal control policy for EBs. The policy aims to jointly optimize energy consumption and traffic efficiency, and is defined by the state space, action space, and reward design described below.

\subsection{State and action}

The state space is constructed from the local observations available to the EB on its current link. We assume that the EB can access the surrounding traffic state, including the distance and speed of the lead vehicle, while the lead-vehicle acceleration is either communicated by a connected lead vehicle or estimated from successive observations. The EB can also receive signal-related information from the traffic controller. The state is defined as
\begin{equation}
    \bm{S}_t\;=\;\left[v_t, a_t, d_l, v_l, a_l, d_{NS}, d_{NI}, \sigma_t, \alpha_t, \delta_t, \mvar{\tau}{rem}, \bm{P}_t \right],\label{eq:state definition}
\end{equation}
where $v_t$ and $a_t$ are the current speed and acceleration of the EB, and $d_l$, $v_l$, and $a_l$ denote the distance to, speed of, and acceleration of the leading vehicle, respectively. The variables $d_{NS}$ and $d_{NI}$ represent the distance to the next stop and the next signalized intersection; $\sigma_t$ indicates whether the EB is in the dwell stage; and $\alpha_t$, as defined in Eq.~\eqref{eq: alpha_j}, indicates whether the bus has departed from the stop and is therefore included in the pressure calculation. The traffic signal state in the EB lane is denoted by $\delta_t$ (green/yellow/red), and $\mvar{\tau}{rem}$ is the residual time of the currently active signal phase.

To encode the pressure information of the downstream intersection, let $\phi_c$ denote the phase serving the EB lane. Its phase pressure is
\begin{equation}
    P_{\phi_c}\;=\;\sum_{(i,o)\in\phi_c}P_{(i,o)}, \label{eq:P-phi-c}
\end{equation}
where $P_{(i,o)}$ is computed by Eq.~\eqref{eq:pressure-3}. Among the remaining phases, we select the one with the largest pressure and denote it by $\phi_m$, with corresponding phase pressure $P_{\phi_m}$. We further define the pressure difference as $\Delta P=P_{\phi_c}-P_{\phi_m}$, and construct the pressure vector
\begin{equation}
    \bm{P}_t\;=\;[P_{\phi_c}, P_{\phi_m}, \Delta P]. \label{eq:state-pressure}
\end{equation}
Together, $d_{NI}$, $\delta_t$, $\mvar{\tau}{rem}$, and $\bm{P}_t$ provide the agent with local signal and priority cues, enabling more informed intersection-approach decisions under adaptive, pressure-driven SPaT. The state design is route-agnostic and can therefore accommodate EBs operating on different routes and intersection layouts. All variables are normalized before being fed to the DRL agent.

The action of the DRL agent is a scalar acceleration command applied to the EB,
\begin{equation}
    A_t=[a^{\mathrm{act}}_t], \qquad
    a^{\mathrm{act}}_t\in[-2.6,2.6]~\text{m/s}^2.
\end{equation}
The action is updated every EB operation step, which is the simulation step, i.e., 1~s in this study.

\subsection{Reward function formulation}
Reward design is pivotal to DRL because it translates high-level control objectives into learnable, step-wise feedback that shapes the policy. 
This agent aims at learning an eco-driving policy that enables the multi-objective and multi-stage operation of EB. Therefore, the reward function $R(\bm{S}_t,a^{\mathrm{act}}_t)$ should be designed accordingly to cover multiple stages and various objectives. According to Eq.~\eqref{eq:objecitve}, the energy consumption and traffic efficiency are the two main objectives of eco-driving.
\subsubsection{Energy consumption reward function}
Energy consumption is a crucial issue of EB, since it directly impacts the travel distance and battery life. 
In order to manage the energy consumption more precisely, an energy consumption model is used to calculate the real-time energy consumed by EB. The model integrates passenger load, road incline, and HVAC system, resulting in a comprehensive mathematical model that calculates the energy consumption based on real-time speed, acceleration, and other external factors. In this study, we employ the model from \cite{yang2024entire} and use it to calculate the instantaneous EB energy consumption.
Based on the real-time operating status of the EB, the step-wise energy reward is defined as
\begin{equation}
R_{\mathrm{eng}}
=
-\frac{[E_t]_+}{E_{\max}},
\qquad
[E_t]_+ = \max(E_t,0),
\label{eq:R-energy}
\end{equation}
where $E_t$ denotes the electrical energy consumed by the EB during control step $t$, and $E_{\max}$ is a normalization constant. Thus, only positive energy consumption is penalized, while regenerative energy recovery does not produce a positive reward.

\subsubsection{Traffic efficiency reward function}
Traffic efficiency is the other major objective of EB operation, which is crucial to reduce traffic congestion and improve public transportation service level. Due to the multi-stage operation feature of EB, the traffic efficiency cannot be evaluated with a unified formulation. Therefore, the efficiency reward function is designed for each operation stage separately. Since efficiency is not relevant when the EB is in the dwell stage, efficiency reward function is not needed for stage 3.  

\noindent\textit{Stage 1: Cruise-stage efficiency and safety reward}

In the cruise stage, the EB travels with the surrounding flow. The reward should encourage steady progress while preserving longitudinal safety under mixed traffic. We therefore combine (i) a small progress incentive to avoid overly conservative cruising, and (ii) a TTC-based soft safety constraint to discourage tailgating.

\textit{a) Progress term.}
Let $v_t$ denote the EB speed at time step $t$. We provide a bounded progress reward only when the vehicle is effectively moving:
\begin{equation}
\label{eq:R-stage1-prog}
R_{\mathrm{prog}}(v_t)
=
\begin{cases}
w_{\mathrm{prog}}\;\min\!\left(1,\dfrac{v_t}{s_{\mathrm{ref}}}\right), & v_t \ge v_{\min}^{\mathrm{prog}},\\[6pt]
0, & v_t < v_{\min}^{\mathrm{prog}},
\end{cases}
\end{equation}
where $w_{\mathrm{prog}}>0$ is the maximum progress reward, $s_{\mathrm{ref}}>0$ is a normalization scale (in m/s), and $v_{\min}^{\mathrm{prog}}>0$ prevents rewarding near-stop crawling.

\textit{b) Safety term.}
Let $\mathrm{TTC}_t$ denote the estimated time-to-collision with the lead vehicle (undefined if no leader exists). Penalties are applied only when $\mathrm{TTC}_t$ falls below preset thresholds $T_{\mathrm{safe}}>T_{\mathrm{hard}}>0$: the \emph{soft band} $[T_{\mathrm{hard}},\,T_{\mathrm{safe}})$ discourages short headways, while the \emph{hard band} $(0,\,T_{\mathrm{hard}})$ imposes stronger penalties to recover safety margins. Specifically,
\begin{equation}
\label{eq:R-stage1-ttc}
\begin{aligned}
R_{\mathrm{ttc}}(\mathrm{TTC}_t)
={}-\,w_{\mathrm{safe}}
&\left(
\frac{[\,T_{\mathrm{safe}}-\mathrm{TTC}_t\,]_+}{T_{\mathrm{safe}}}
\right)^{p} \\
&-\,w_{\mathrm{hard}}
\left(
\frac{[\,T_{\mathrm{hard}}-\mathrm{TTC}_t\,]_+}{T_{\mathrm{hard}}}
\right)^{p},
\end{aligned}
\end{equation}
where $w_{\mathrm{safe}}>0$ and $w_{\mathrm{hard}}>0$ weight the two bands, and $p\in\{1,2\}$ controls the curvature (linear or quadratic). For numerical robustness, we set $R_{\mathrm{ttc}}=0$ when $\mathrm{TTC}_t$ is undefined, non-positive, or exceeds a large cutoff (e.g., $\mathrm{TTC}_t \ge 10$~s), in which case collision risk is negligible.

\textit{c) Stage-1 reward.}
The cruise-stage reward is then
\begin{equation}
\label{eq:R-stage1}
R_{\mathrm{stage1}} \;=\; R_{\mathrm{prog}}(v_t) \;+\; R_{\mathrm{ttc}}(\mathrm{TTC}_t),
\end{equation}
which promotes forward motion while softly enforcing safe car-following.

\noindent\textit{Stage 2: Stop-approach stage efficiency reward}

During the stop-approach stage, the agent must ensure a gentle, low-speed entry into the bus stop. Similar to the reward design in \cite{zhang2024eco}, we predict the entry speed at the stop line using constant-acceleration kinematics over the remaining approach distance. With the current speed $v_t$ and acceleration $a_t$, the distance to the next stop $d_{NS}$, and a look-ahead cap $d_L$, we define the effective approach distance
\begin{equation}
\mvar{d}{stop} \;=\; \min\{\,d_{NS},\, d_L\,\},
\end{equation}
and the predicted entry speed
\begin{equation}
\label{eq:v_entry}
v_{\mathrm{entry}} \;=\; \sqrt{\,v_t^2 + 2\,a_t \,\mvar{d}{stop}\,},
\end{equation}
where the discriminant $v_t^2+2a_t\mvar{d}{stop}$ must be nonnegative. Let $V_{\mathrm{safe}}$ denote the maximum acceptable entry speed in the stop area (e.g., $V_{\mathrm{safe}}=2.6\,\mathrm{m/s}$). The reward is then defined as a piecewise shaping term:
\begin{equation}
\label{eq:R-stage2}
\mvar{R}{stage2} \;=\;
\begin{cases}
-1, & v_t^2 + 2 a_t \mvar{d}{stop} < 0 \quad (\text{unsafe}),\\[4pt]
\displaystyle 2\,\frac{V_{\mathrm{safe}} - v_{\mathrm{entry}}}{V_{\mathrm{safe}}}, & 0 \le v_{\mathrm{entry}} \le V_{\mathrm{safe}},\\[10pt]
-3, & v_{\mathrm{entry}} > V_{\mathrm{safe}}.
\end{cases}
\end{equation}
Thus the agent receives a positive linear reward for planning a sufficiently slow entry (maximal at $v_{\mathrm{entry}}{=}0$ and tapering to $0$ at $V_{\mathrm{safe}}$), a strong penalty if the predicted entry speed exceeds $V_{\mathrm{safe}}$, and a mild penalty when the kinematic prediction is invalid (negative discriminant). This shaping promotes safe, smooth docking while preserving learnability.

\noindent\textit{Stage 4: Intersection-crossing stage efficiency reward}

Intersection crossing is the most delicate stage because the SPaT observed by the EB is time-varying and stochastic under Priority-MP control. Rather than using a future phase table, we rely only on the current residual time and local pressure cues. The idea is to estimate the likelihood that the current green will be extended, convert it into an expected usable green window, and compare that window with the EB’s expected time of arrival (ETA). The resulting reward consists of a guidance term and a sparse event bonus for passing on green.

\textit{a) Time-aware green-extension probability.}
Recall the pressure state vector $\bm{P}_t=[P_{\phi_c}, P_{\phi_m}, \Delta P]$ from Eqs.~\eqref{eq:P-phi-c}--\eqref{eq:state-pressure}. We define the normalized phase advantage
\begin{equation}
\label{eq:adv}
\mathrm{adv} \;=\; \frac{\Delta P}{|P_{\phi_c}|+|P_{\phi_m}|+\varepsilon}\in(-1,1),
\end{equation}
which measures how strongly the current phase $\phi_c$ is favored over its main competitor $\phi_m$ (normalization equalizes scale across intersections and $\varepsilon>0$ prevents division by zero). To reflect decision urgency, we introduce the time-aware weight
\begin{equation}
\label{eq:time-weight}
\mvar{w}{adv} \;=\; \max\!\Big\{\,w_{\min},\; 1-\min\!\big(1,\;\mvar{\tau}{rem}/T_c\big)\,\Big\}\in[w_{\min},1],
\end{equation}
where $\mvar{\tau}{rem}$ is the residual phase time and $T_c$ is the signal control step. A pressure-based extension likelihood is then obtained via
\begin{equation}
\label{eq:p_pressure}
p_{\mathrm{press}} \;=\; \frac{1}{1+e^{-k\,\mathrm{adv}}}\in(0,1),
\end{equation}
a smooth, bounded transformation of $\mathrm{adv}$ whose slope $k>0$ controls sensitivity around $\mathrm{adv}=0$. Finally, we blended it with an uninformed prior using the time-aware weight:
\begin{equation}
\label{eq:p_ext}
p_{\mathrm{ext}} \;=\; (1-\mvar{w}{adv})\cdot 0.5 \;+\; \mvar{w}{adv}\cdot p_{\mathrm{press}} \;\in\;[0,1].
\end{equation}
Thus, early in the phase the estimate stays conservative, while near the switching instant it relies more on pressure information.

\textit{b) Expected extra green.}
Treating $p_{\mathrm{ext}}$ as the per-step probability that the current phase continues for one more control step of length $T_c$, we estimate the expected extra green beyond $\mvar{\tau}{rem}$ as
\begin{equation}
\label{eq:exp-green}
\mathbb{E}[\Delta G] \;=\; T_c \sum_{i=1}^{H}\rho^{\,i-1}\,p_{\mathrm{ext}}^{\,i}
\;=\; T_c\;\frac{p_{\mathrm{ext}}\big(1-(\rho\,p_{\mathrm{ext}})^H\big)}{1-\rho\,p_{\mathrm{ext}}},
\end{equation}
where $\rho\in(0,1)$ discounts long extension chains and $H$ truncates the horizon.

\textit{c) Guidance reward using ETA vs.\ expected green.}
We compare the vehicle’s arrival time with the usable green window implied by the current SPaT. Let $v_t$ be the current speed, $d_{NI}$ the distance to the stop line, and
\begin{equation}
\label{eq:eta}
\mathrm{ETA} \;=\; \frac{d_{NI}}{\max(v_t,\epsilon)} ,\quad \epsilon>0.
\end{equation}
Under a green signal ($\mathrm{G}$), the expected usable green is
\begin{equation}
\label{eq:Gexp-green}
G_{\mathrm{exp}}^{(\mathrm{G})} \;=\; [\mvar{\tau}{rem}-\mathrm{marg}]_++\mathbb{E}[\Delta G],
\end{equation}
where $\mathrm{marg}>0$ is a small scheduling margin. Let $\mathrm{slack}=G_{\mathrm{exp}}^{(\mathrm{G})}-\mathrm{ETA}$ that measures catchability, and the corresponding shaping reward is
\begin{equation}
\label{eq:R-green}
R_{\mathrm{guide}}^{(\mathrm{G})} \;=\;
\begin{cases}
\displaystyle \mvar{r}{guide}^{+}\cdot\min\!\left(1,\;\frac{\mathrm{slack}+\mathrm{marg}}{G_{\mathrm{exp}}^{(\mathrm{G})}+\mathrm{marg}+\epsilon}\right), & \mathrm{slack}\ge 0,\\[8pt]
\displaystyle \mvar{r}{guide}^{-}\cdot\min\!\left(1,\;\frac{-\mathrm{slack}}{\mathrm{marg}+\mathrm{ETA}}\right), & \mathrm{slack}<0,
\end{cases}
\end{equation}
so that feasible catch-the-green arrivals are rewarded (up to $\mvar{r}{guide}^{+}$), while late arrivals are gently penalized (down to $\mvar{r}{guide}^{-}$); the $\min(\cdot,1)$ caps magnitude to stabilize learning.

For yellow ($\mathrm{Y}$) or red ($\mathrm{R}$), we encourage planning beyond the imminent switch rather than rushing. We therefore define
\begin{equation}
\label{eq:trem-adj}
\widetilde{t}_{\mathrm{rem}} \;=\;
\begin{cases}
\mvar{\tau}{rem} + T_c, & \mathrm{Y},\\
\mvar{\tau}{rem} + T_c\,(1-p_{\mathrm{ext}}), & \mathrm{R},
\end{cases}
\quad
\mvar{e}{arr} \;=\; \mathrm{ETA} - (\widetilde{t}_{\mathrm{rem}}+\mathrm{marg}),
\end{equation}
and
\begin{equation}
\label{eq:R-redyellow}
R_{\mathrm{guide}}^{(\mathrm{Y/R})} \;=\;
\begin{cases}
\displaystyle \mvar{r}{guide}^{+}\cdot\min\!\left(1,\;\frac{\mvar{e}{arr}}{\mvar{e}{arr}+\mathrm{marg}}\right), & \mvar{e}{arr}\ge 0,\\[8pt]
\displaystyle \mvar{r}{guide}^{-}\cdot\min\!\left(1,\;\frac{-\mvar{e}{arr}}{\mathrm{marg}+\widetilde{t}_{\mathrm{rem}}+\epsilon}\right), & \mvar{e}{arr}<0.
\end{cases}
\end{equation}
This construction rewards feasible green passage under $\mathrm{G}$ and discourages risky rushing under $\mathrm{Y}/\mathrm{R}$.

\textit{d) Event bonus and total reward.}
Let $\mathbb{I}_{\mathrm{passG}}\in\{0,1\}$ indicate whether the EB crosses the stop line under green at the current step. We add the sparse event bonus
\begin{equation}
\label{eq:R-event}
R_{\mathrm{event}} \;=\; w_{\mathrm{pass}}\;\mathbb{I}_{\mathrm{passG}},
\end{equation}
and define the stage-4 reward as
\begin{equation}
\label{eq:R-stage4}
R_{\mathrm{stage4}} \;=\;
\begin{cases}
R_{\mathrm{guide}}^{(\mathrm{G})} + R_{\mathrm{event}}, & \text{if the signal state is } \mathrm{G},\\[3pt]
R_{\mathrm{guide}}^{(\mathrm{Y/R})} + R_{\mathrm{event}}, & \text{if the signal state is } \mathrm{Y}\text{ or }\mathrm{R}.
\end{cases}
\end{equation}

Overall, the stage-4 reward follows a guidance-plus-event structure: the guidance term steers ETA toward a feasible green window inferred from local pressure cues and residual time, while the event term reinforces successful green passage. Because the design depends only on local SPaT information and pressures, it remains lightweight and transferable across intersections.

\noindent\textit{Overall efficiency reward}

With the stage-wise rewards defined above, the agent must first identify the current operating stage and then apply the corresponding reward. We use a lightweight rule-based stage classifier that gives priority to \emph{stop-approach} when a stop is imminent, switches to \emph{intersection-crossing} when the signal becomes the nearer target and lies within a speed-dependent preparation distance, and otherwise keeps the EB in \emph{cruise}. Stage identification is performed only when the EB is not dwelling at a stop, since the dwell stage is handled separately.

Given the current speed $v_t$, the preparation distance for signal approach is defined as
\begin{equation}
\label{eq:tls-thr-simple}
\begin{aligned}
D_{\mathrm{prep}}(v_t)
=\,&\mathrm{clip}\!\Bigl(
\max\!\Bigl\{
v_t t_{\mathrm{win}}, \\
&
v_t \mvar{\tau}{reac}
+\frac{v_t^{2}}{2a_{\mathrm{comf}}}
+\mvar{d}{static}
+\text{buffer}
\Bigr\},
[\underline{d}_{\mathrm{isec}},\,\bar{d}_{\mathrm{isec}}]
\Bigr).
\end{aligned}
\end{equation}
Here, $t_{\mathrm{win}}$ is a short preview window, $\mvar{\tau}{reac}$ is the reaction time, $a_{\mathrm{comf}}$ is a comfortable deceleration, and $\mvar{d}{static}$ is a static safety distance with an additional buffer. This construction combines a time-based look-ahead with a physics-based stopping distance, so that the EB enters the signal-approach stage early enough either to align with a feasible green or to brake comfortably if not.

Let $d_{NI}\in(0,+\infty]$ denote the distance to the next signalized intersection (set to $+\infty$ if none exists on the current segment), and let $d_{NS}\in(0,+\infty]$ denote the distance to the next stop (also set to $+\infty$ if none exists). With a near-stop threshold $d_L$, the stage is determined as
\begin{equation}
\label{eq:stage-map-simple}
\mathrm{stage} =
\begin{cases}
2, & d_{NS}\le d_{L}\ , \ d_{NS}<d_{NI},\\[3pt]
4, & d_{NI}\le D_{\mathrm{prep}}(v_t)\ ,\ d_{NI}<d_{NS},\\[3pt]
1, & \text{otherwise}.
\end{cases}
\end{equation}
This rule ensures that stop docking takes precedence when necessary, while intersection-related decision making is activated only when the signal is the nearer target and the EB has entered the preparation zone.

Let $q_t\in\{1,2,4\}$ denote the current stage. To discourage overly conservative behavior and reflect schedule-reliability requirements, we include a small per-step time cost $r_{\mathrm{time}}<0$ in the efficiency reward. The overall efficiency reward is then
\begin{equation}
\label{eq:overall-eff-cases}
\mvar{R}{eff} =
\begin{cases}
R_{\mathrm{stage}1} + r_{\mathrm{time}}, & q_t=1,\\
R_{\mathrm{stage}2} + r_{\mathrm{time}}, & q_t=2,\\
R_{\mathrm{stage}4} + r_{\mathrm{time}}, & q_t=4,
\end{cases}
\end{equation}
where $R_{\mathrm{stage}1}$, $R_{\mathrm{stage}2}$, and $R_{\mathrm{stage}4}$ are defined in
\eqref{eq:R-stage1}, \eqref{eq:R-stage2}, and \eqref{eq:R-stage4}, respectively.

\subsubsection{Reward functions for other objectives}
Beyond the stage-conditioned \emph{efficiency} terms, we include auxiliary rewards to promote ride comfort, safety, and operational regularity. These terms are deliberately simple, bounded, and additive so they can be tuned without destabilizing learning.

\paragraph{Jerk smoothing (comfort)}
Let $a_t$ and $a_{t-1}$ denote longitudinal accelerations and $\Delta t$ the EB control step. The jerk is
\begin{equation}
j_t = \frac{a_t-a_{t-1}}{\Delta t},
\end{equation}
and we penalize excessive jerk quadratically with a normalization constant $J_0>0$:
\begin{equation}
\label{eq:r-jerk}
\mvar{R}{jerk} = -\,\frac{j_t^{\,2}}{J_0}\;\;\le 0.
\end{equation}
In implementation, we use $\Delta t=1.0$ s and $J_0=27.04$ (so that $\sqrt{J_0}\approx \frac{2.6-(-2.6)}{1}=5.2~\mathrm{m/s^3}$), which sets the scale of this penalty term under the adopted action bounds.

\paragraph{Collision penalty (safety)}
Let $\mathrm{col}_t\!\in\!\{0,1\}$ indicate whether a collision is detected at time $t$. A hard safety penalty is applied:
\begin{equation}
\label{eq:r-collision}
\mvar{R}{col} = -\,w_{\mathrm{col}}\,\mathrm{col}_t,\qquad w_{\mathrm{col}}>0,
\end{equation}
with $w_{\mathrm{col}}{=}10$ by default. If no collision occurs, the term is zero.

\paragraph{Unscheduled stop penalty (operations)}
To discourage full stops outside designated berths, we add a small penalty when the vehicle is (nearly) stationary away from a bus stop. With the current speed $v_t$, the distance to the next bus stop $d_{NS}$, and a small threshold $\varepsilon_v$, we have:
\begin{equation}
\label{eq:r-spurious-stop}
\mvar{R}{usp} =
\begin{cases}
-\,w_{\mathrm{usp}}, & \text{if } d_{NS}>0 \ \text{and}\ v_t \le \varepsilon_v,\\
0, & \text{otherwise},
\end{cases}
\end{equation}
with $w_{\mathrm{usp}}{=}5.0$ and $\varepsilon_v{=}0.5$ m/s by default. Note that the penalty is not applied when the EB is stuck in congestion or waiting for the traffic light.

\paragraph{Harsh braking penalty (comfort/safety)}
Finally, overly aggressive braking receives a small penalty:
\begin{equation}
\label{eq:r-brake}
\mvar{R}{hb} =
\begin{cases}
-\,w_{\mathrm{hb}}, & \text{if } a_t < -a_{\mathrm{thr}},\\
0, & \text{otherwise},
\end{cases}
\end{equation}
where $a_{\mathrm{thr}}{=}2.0$\,m/s$^2$ and $w_{\mathrm{hb}}{=}0.5$ by default.

\subsubsection{Overall reward}
The auxiliary (non-efficiency) component is the sum of the above terms,
\begin{equation}
\label{eq:r-other}
R_{\mathrm{oth}} \;=\; R_{\mathrm{jerk}} \;+\; R_{\mathrm{col}} \;+\; R_{\mathrm{usp}} \;+\; R_{\mathrm{hb}} \;\;\le 0,
\end{equation}
and the overall training signal combines the stage-conditioned efficiency reward with energy and comfort/safety weights:
\begin{equation}
\label{eq:r-total}
R_t \;=\; \omega_e\, R_{\mathrm{eng}} \;+\; \omega_t\, \mvar{R}{eff} \;+\; \omega_o\, \mvar{R}{oth},
\end{equation}
where $(\omega_e,\omega_t,\omega_o)$ and the constants in \eqref{eq:r-jerk}–\eqref{eq:r-brake} are tunable parameters. This design keeps safety terms always active, enforces smoothness through jerk shaping, and discourages operationally undesirable full stops outside berths, while leaving most of the long-horizon timing behavior to the efficiency rewards.

\subsection{CTDE with parameter-sharing training framework}

Eco-driving requires continuous control under mixed continuous--discrete observations and dynamically changing traffic conditions. We therefore adopt Soft Actor--Critic (SAC) \cite{haarnoja2018soft}, an off-policy maximum-entropy algorithm well suited to continuous action spaces. SAC combines stochastic policy learning, replay-based sample reuse, and double-Q value estimation with target networks, which together provide strong stability and sample efficiency in nonstationary environments. These properties make SAC an appropriate optimization backbone for the proposed pressure-aware EB eco-driving problem.

To further improve data efficiency and transferability, we employ a centralized-training, decentralized-execution (CTDE) framework with parameter sharing. Instead of training a separate policy for each route or vehicle, we train a single shared SAC actor--critic using experiences collected from multiple buses. This is enabled by the route-agnostic observation design, which encodes local SPaT and stage information so that the same policy can be applied across different routes and layouts.

During training, the ego vehicle is switched every period $\Delta\mvar{T}{sw}$ (or when the current ego EB leaves the network), after which interaction continues on the newly assigned bus. This mechanism exposes the learner to diverse operating conditions, including different link geometries, stop spacing, traffic demand, and pressure patterns, thereby enriching the replay buffer while keeping the controller architecture unchanged. The shared actor and critics are then updated off-policy from the aggregated replay memory.
At deployment, execution is fully decentralized: the trained policy runs independently on each EB using only local observations, without requiring online coordination among buses. This CTDE with parameter-sharing framework therefore improves training efficiency and generalization while preserving real-time onboard implementability. The pseudocode of the training procedure is given in Algorithm~\ref{alg:prieco-sac}. Although SAC is adopted in this study, the overall framework is not restricted to SAC and can be extended to other DRL algorithms.

\begin{algorithm}[!h]
\small
\caption{PriEco\textendash DRL training with SAC under CTDE and parameter sharing}
\label{alg:prieco-sac}
\begin{algorithmic}[1]  
\Require simulation interval length $T$, ego\hyp switch period $\Delta T_{\mathrm{sw}}$, signal control cycle $\mvar{T}{c}$, replay capacity $|\mathcal{R}|$
\Ensure shared actor $\pi_\theta$, shared double critics $Q_{\psi_1}, Q_{\psi_2}$
\State Initialize the actor network $\pi_\theta(a^{\mathrm{act}}|\bm{S})$, critic networks $Q_{\psi_1},Q_{\psi_2}$, target networks $\bar Q_{\psi_1},\bar Q_{\psi_2}$
\State Initialize shared replay buffer $\mathcal{R}$, entropy temperature $\mvar{\alpha}{ent}$
\For{simulation run $=1,2,\dots$}
  \State  Reset the environment; set Priority-MP  control by \eqref{eq:pressure-3}-\eqref{eq:priority_MP_with_ts}
  \State Identify EB set $\mathcal{B}$ on the network; sample ego $b \sim \mathrm{Unif}(\mathcal{B})$; obtain initial state $\bm{S}_1$
  \State Set the initial switch time $t_{\mathrm{next}} \gets \Delta T_{\mathrm{sw}}$
  \For{$t=1$ to $T$}
    \State Sample an action from the policy $a^{\mathrm{act}}_t \sim \pi_\theta(\cdot|\bm{S}_t)$ 
    \State Project $a^{\mathrm{act}}_t$ by rule layer to ensure  feasible acceleration, speed, and headway
    \State Perform EB eco-driving strategy by applying acceleration action $a^{\mathrm{act}}_t$ on ego $b$ 
    \State Simulate the environment by one step
    \State Calculate the movement pressures $P_{(i,o)}$ using \eqref{eq:pressure-3}
    \If{$t \bmod T_c = 0$} 
    \State Perform Priority-MP signal control \eqref{eq:priority_MP_with_ts} determining the phases for the next control cycle
    \EndIf
    \State Get EB status, traffic environment, and SPaT information $\delta_t,\mvar{\tau}{rem},\bm{P}_t$; construct state  $\bm{S}_{t+1}$
    \State Compute energy consumption reward $R_{\mathrm{eng}}$ by \eqref{eq:R-energy}
    \State Compute $D_{\mathrm{prep}}(v_t)$ by \eqref{eq:tls-thr-simple}
    \State Determine stage $q_t \in \{1,2,4\}$ by \eqref{eq:stage-map-simple}
    \State Set the efficiency reward $\mvar{R}{eff}$ according to \eqref{eq:overall-eff-cases}
    \State Compute $\mvar{R}{jerk}$ \eqref{eq:r-jerk}, $\mvar{R}{col}$ \eqref{eq:r-collision}, $\mvar{R}{usp}$ \eqref{eq:r-spurious-stop}, $\mvar{R}{hb}$ \eqref{eq:r-brake}
    \State Calculate overall reward by \eqref{eq:r-total} 
    \State Push $(\bm{S}_t,a^{\mathrm{act}}_t,R_t,\bm{S}_{t+1})$ into $\mathcal{R}$ (evict oldest if full)
    \If{update SAC model}
    \For{$k=1$ to $K$} 
      \State Sample mini\hyp batch $\{(\bm{S},a^{\mathrm{act}},R,\bm{S}')\}$ from $\mathcal{R}$
      \State $\tilde a' \sim \pi_\theta(\cdot|\bm{S}')$;
      \Statex\qquad\qquad\qquad $y \gets r+\gamma\!\left(\min_{j=1,2}\bar Q_{\psi_j}(\bm{S}',\tilde a')-\mvar{\alpha}{ent}\log\pi_\theta(\tilde a'|\bm{S}')\right)$
      \State Update critics: 
      \Statex\qquad\qquad\qquad minimize $\sum (Q_{\psi_j}(\bm{S},a)-y)^2$ for $j=1,2$
      \State $\tilde a \sim \pi_\theta(\cdot|\bm{S})$;\quad
             update actor: 
        \Statex\qquad\qquad\qquad maximize $\sum \left(\min_{j} Q_{\psi_j}(\bm{S},\tilde a)-\mvar{\alpha}{ent}\log\pi_\theta(\tilde a|\bm{S})\right)$
      \State Auto\hyp tune $\mvar{\alpha}{ent}$; soft\hyp update 
      \Statex\qquad\qquad\qquad$\bar Q_{\psi_j}\gets \tau Q_{\psi_j}+(1-\tau)\bar Q_{\psi_j}$ for $j=1,2$
    \EndFor
    \EndIf
    \If{$t \ge t_{\mathrm{next}}$ \textbf{or} ego $b$ is outside the network } 
      \State Update EB set $\mathcal{B}$ on the network
      \State Sample new ego $b \sim \mathrm{Unif}(\mathcal{B})$;\quad $t_{\mathrm{next}} \gets t + \Delta T_{\mathrm{sw}}$
    \EndIf
    \State $\bm{S}_t \gets \bm{S}_{t+1}$
  \EndFor
\EndFor
\end{algorithmic}
\end{algorithm}


\section{Numerical experiments} \label{sec:case-study}
\subsection{Experimental settings}
\begin{figure}[h]
    \centering
    \includegraphics[width=1.0\linewidth]{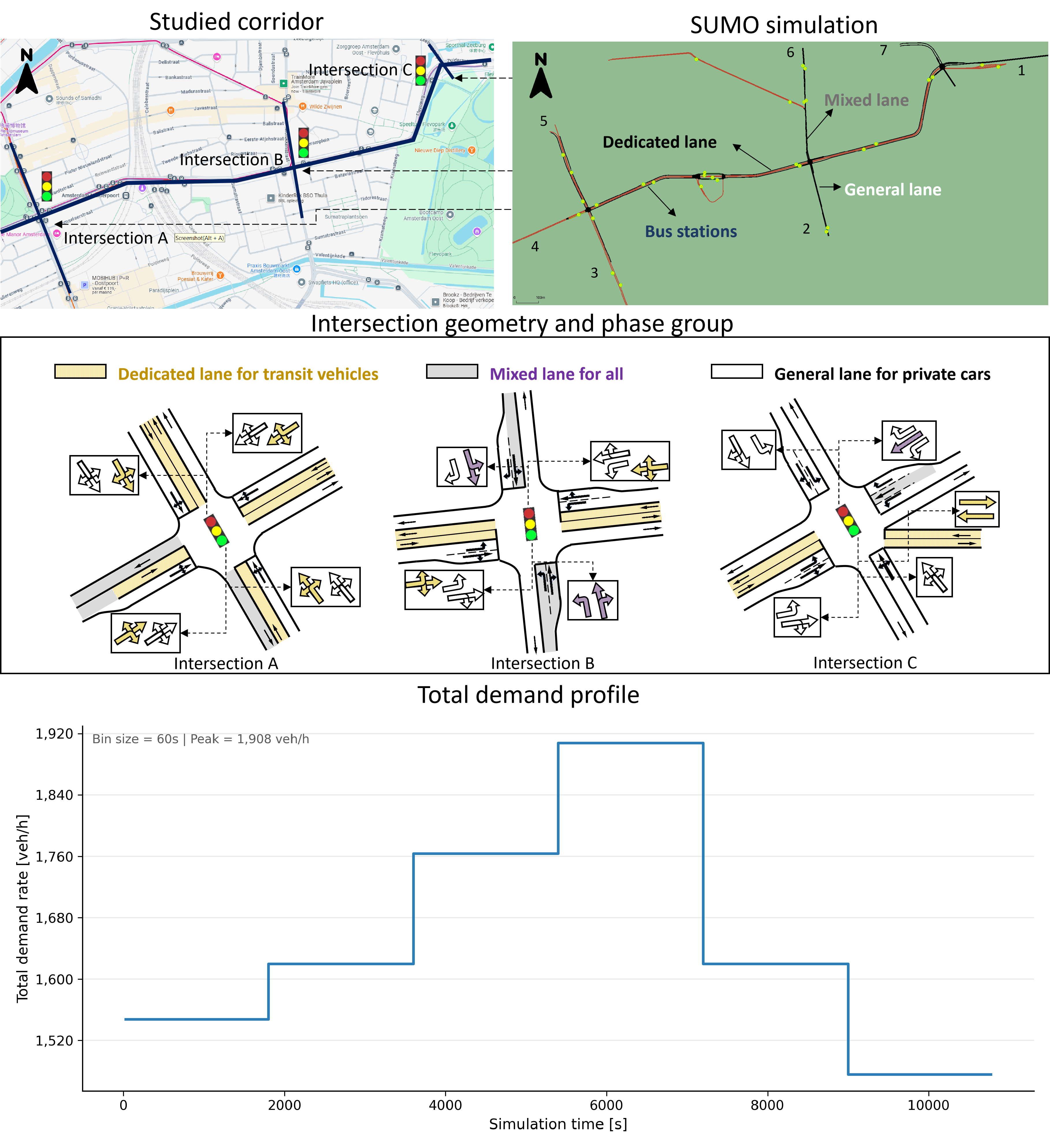}
    \vspace{-0.5cm}
    \caption{Layout of the studied corridor in Amsterdam with three intersections and the total demand in the network}
    \label{fig:Case-study}
\end{figure}

We evaluate the proposed PriEco-DRL framework on a real-world corridor in Amsterdam. The corridor contains three signalized intersections, eight bus lines, and 18 bus stops, as shown in Fig.~\ref{fig:Case-study}. The corresponding SUMO simulation environment includes dedicated lanes for buses and trams, mixed-use lanes, and general-purpose lanes for private vehicles. Following the real-world signal layout, movements from the same incoming link are grouped into one phase group. EBs on dedicated lanes are granted higher priority during green, which in practice can be implemented by warning private cars before transit passage.

\begin{figure}[h]
    \centering
    \includegraphics[width=1.0\linewidth, trim={0cm, 3.5cm, 0cm, 3.5cm}, clip]{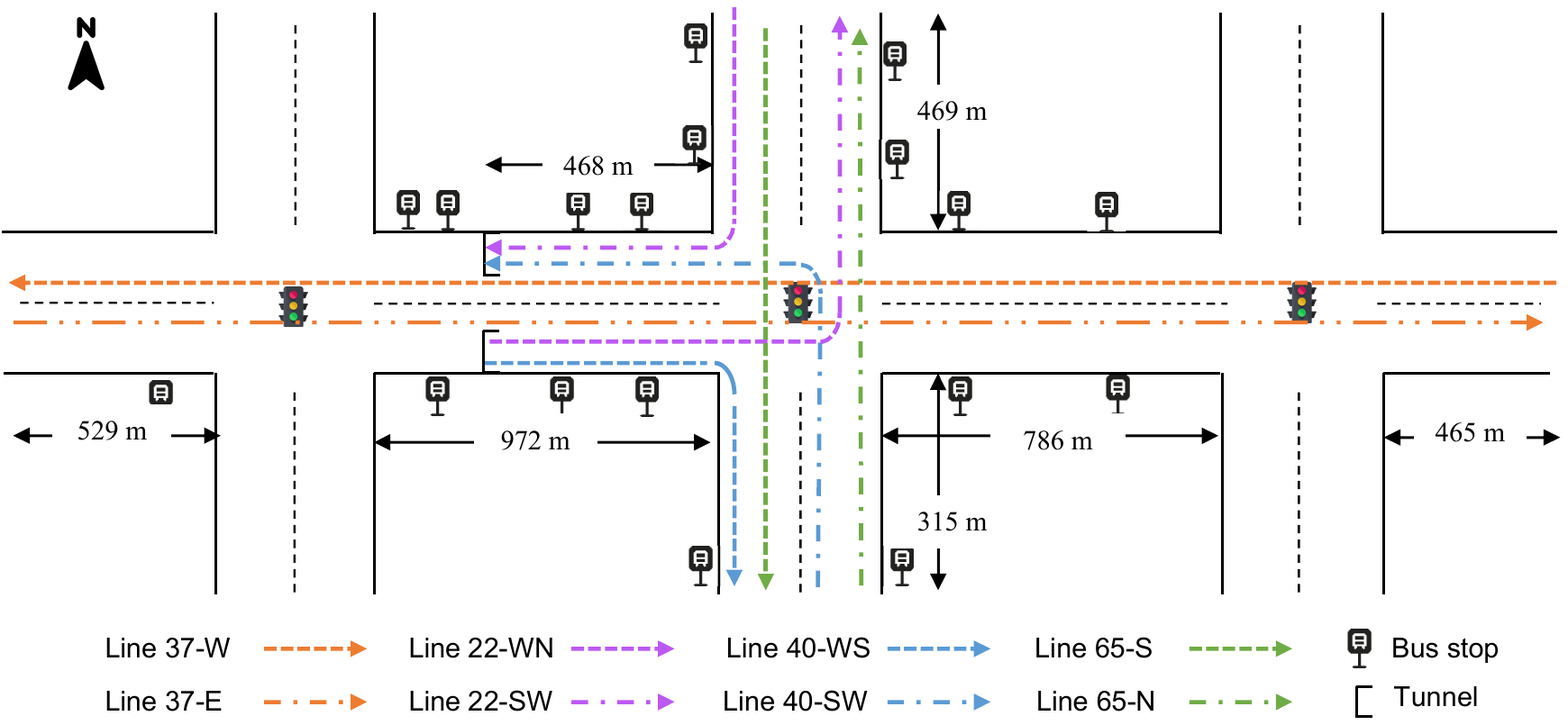}
    \vspace{-0.5cm}
    \caption{Bus routes and stops distribution. Route labels use line number plus direction abbreviation.}
    \label{fig:bus-route}
\end{figure}

The corridor has eight boundary links, resulting in eight origins/destinations and 56 OD pairs for private-car traffic. Over the 3-hour simulation horizon, the demand on the two main-road OD pairs (1--5 and 5--1) increases from 72 veh/h to 252 veh/h and then decreases to 36 veh/h, while all other OD pairs are fixed at 30 veh/h. The total corridor demand therefore varies between approximately 1500 and 2000 veh/h.
Public-transit departures follow the real-world headways, i.e., 10 or 20 min depending on the line. Figure~\ref{fig:bus-route} shows the route geometry and stop distribution, together with route lengths and stop counts. For readability, the bus lines are denoted by line number and direction: Line 37-W/37-E, Line 40-WS/40-SW, Line 22-WN/22-SW, and Line 65-S/65-N, where W, E, N, and S indicate west, east, north, and south, respectively. EBs are dispatched evenly over time within each run. During the 3-hour simulation, 12 trips are completed on Lines 37-W, 37-E, 65-S, and 65-N; 9 trips on Lines 40-WS and 40-SW; and 18 trips on Lines 22-WN and 22-SW.

For signal control, we consider three methods: (i) fixed-time control, with deterministic phase durations; (ii) actuated control in SUMO, which adjusts green times using detector inputs; and (iii) the proposed Priority-MP controller described in Section~\ref{sec:MP}. For EB longitudinal control, we consider: (i) the SUMO IDM model \cite{treiber2000congested}; (ii) a rule-based eco-driving module, denoted as G2, that combines GLOSA \cite{barth2011dynamic} and GLODTA; and (iii) the proposed PriEco-DRL policy.

To construct coupled baselines, we integrate G2 with each signal controller in a plug-and-play manner. At each step, the signal controller updates the current indication, while each bus computes a desired speed from local SPaT information using GLOSA; at stops, GLODTA may add a holding time (subject to minimum dwell constraints) to improve downstream green alignment. This yields the following comparison methods:
\begin{itemize}
    \item FT-IDM: fixed-time signal control with IDM;
    \item FT-G2: fixed-time signal control with GLOSA--GLODTA;
    \item AC-IDM: actuated control with IDM;
    \item AC-G2: actuated control with GLOSA--GLODTA;
    \item PMP-IDM: Priority-MP control with IDM;
    \item PMP-G2: Priority-MP control with GLOSA--GLODTA;
    \item PriEco-DRL: the proposed integrated framework.
\end{itemize}

We evaluate the methods from both the network and EB perspectives. Network-level metrics include: (i) average vehicle delay (s), (ii) network vehicle count, (iii) maximum queue count, and (iv) maximum spillover count. EB-level metrics include: (i) EB delay (s), defined as time loss relative to free-flow travel time over all completed trips; (ii) EB energy consumption (kWh), defined as the total electrical energy consumed by all completed EB trips, including traction and auxiliary loads; and (iii) EB unscheduled stops, defined as the number of full stops outside designated bus stops.
The SAC hyperparameters are listed in Table~\ref{tab:sac-hparams}, and
all the key simulation parameter values are presented in Table~\ref{tab:key_params}.
To examine the energy--travel-time trade-off, we test three reward-weight settings: PriEco-DRL-T with $(\omega_e,\omega_t,\omega_o)=(1,2,1)$, PriEco-DRL-E with $(2,1,1)$, and PriEco-DRL-B with $(1,1,1)$.
\begin{table}[t]
\centering
\caption{SAC training hyperparameters used in this study}
\label{tab:sac-hparams}
\footnotesize
\begin{tabular}{@{}l l@{}}
\toprule
\textbf{Category} & \textbf{Setting (SAC)} \\
\midrule
\multicolumn{2}{l}{\textit{Training \& Optimization}} \\
Total timesteps & $5 \times 10^{5}$ \\
Learning rate  & $3\times10^{-4}$ \\
Discount factor  & $0.99$ \\
Batch size & 512 \\
Replay buffer size & $200{,}000$ \\
Train frequency / Gradient steps & 1 / 1 \\
Soft update coefficient & 0.002 \\
Entropy coefficient & \texttt{"auto"} (automatic temperature) \\
Target entropy & \texttt{"auto"} $\approx -\mathrm{dim}(A_t)$ \\
\midrule
\multicolumn{2}{l}{\textit{Policy / Networks}} \\
Policy type & \texttt{MlpPolicy} \\
Hidden layers & [256, 256] \\
Activation & ReLU\\
\bottomrule
\end{tabular}

\end{table}

\begin{table}[htbp]
\centering
\caption{Key simulation and reward parameters}
\label{tab:key_params}
\setlength{\tabcolsep}{6pt}
\begin{tabular}{ccccccc}
\toprule
$T_c$ & $\Delta T_{\mathrm{sw}}$ & $\mvar{w}{prog}$ & $\mvar{s}{ref}$ & $\mvar{v}{min}^{\mathrm{prog}}$ & $\mvar{w}{safe}$ & $\mvar{w}{hard}$ \\
10 s & 500 s & 1.0 & 10 m/s & 0.5 & 1.0 & 3.0 \\
\midrule
$\mvar{T}{safe}$& $\mvar{T}{hard}$ & $d_L$ & $\mvar{r}{guide}^+$ & $\mvar{r}{guide}^-$ & $\mvar{w}{pass}$ & $a_{\mathrm{comf}}$ \\
3.5 s & 1.5 s & 50 m & 3.0 & -3.0 & 5.0 & 2.0 m/s$^2$ \\
\midrule
$\mvar{w}{min}$ & $\underline{d}_{\mathrm{isec}}$ & $\bar{d}_{\mathrm{isec}}$ & $\mvar{t}{win}$ & $\mvar{\tau}{reac}$ & $\mvar{d}{static}$ & $\mvar{r}{time}$\\
0.15 & 30 m & 120 m & 10 s & 2.0 s & 5.0 m & -0.1\\
\bottomrule
\end{tabular}
\end{table}

\subsection{Model training}
\begin{figure}
    \centering
    \includegraphics[width=1.0\linewidth, trim={0cm, 0cm, 0cm, 0cm}, clip]{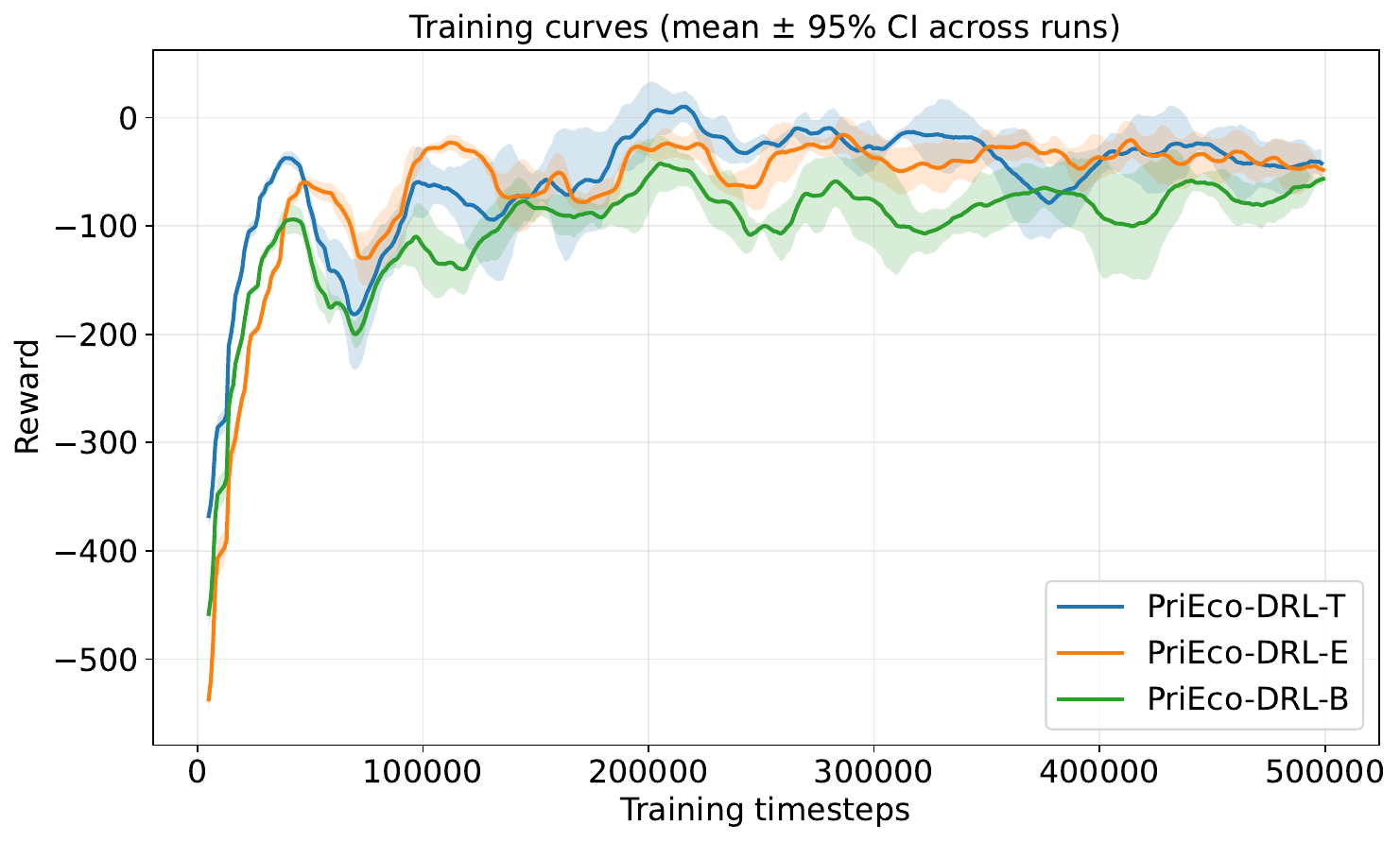}
    \vspace{-0.5cm}
    \caption{Training curves of PriEco-DRL variants under SAC with a fixed ego-switch interval $\Delta T_{\mathrm{sw}}=500$ s. The Y-axis shows the episode return (episode reward), where each episode corresponds to one $500$-s simulation segment between ego switches; the X-axis is the number of environment timesteps. Solid lines report the mean over 5 random seeds, and shaded bands indicate the 95\% confidence interval.} 
    \label{fig:reward}
\end{figure}

During each training run, we optimize a SAC agent under a CTDE, parameter-sharing setup for a total of $5\times10^{5}$ environment steps.
Only a subset of bus lines are used as ego vehicles during training, namely Line 37-W, Line 22-WN, Line 40-SW, and Line 65-S. 
The training process is organized into fixed-length episodes with a constant interval $\Delta T_{\mathrm{sw}}=500$ s.
In each episode, the ego vehicle remains active for 500 seconds before switching to another ego bus. This protocol allows for dynamic switching of buses while ensuring consistent interaction data across multiple routes.
At test stage, the learned policy is deployed to any EB across the full network using only local observations. 
This protocol explicitly separates the training and evaluation distributions by training on a subset of routes and evaluating on the entire network, thereby assessing transferability beyond a single-route setting with concurrent buses from different lines.
For each PriEco-DRL variant (PriEco-DRL-T/E/B), training is repeated five times with distinct random seeds. 
Figure~\ref{fig:reward} reports the learning curves as the mean episode return across runs with a shaded 95\% confidence interval. 
Since the episode interval is fixed, episode returns are directly comparable across runs and variants; we additionally smooth the curves for visualization while preserving the underlying trend. 
As shown, all variants exhibit rapid performance improvement in early training and then gradually stabilize, indicating that SAC converges to a steady policy under the proposed training setup.
The hyperparameters of SAC are listed in Table~\ref{tab:sac-hparams}.

\subsection{Results and analysis}

\subsubsection{Overall performance}
\begin{figure}[]
    \centering
    \includegraphics[width=1.0\linewidth]{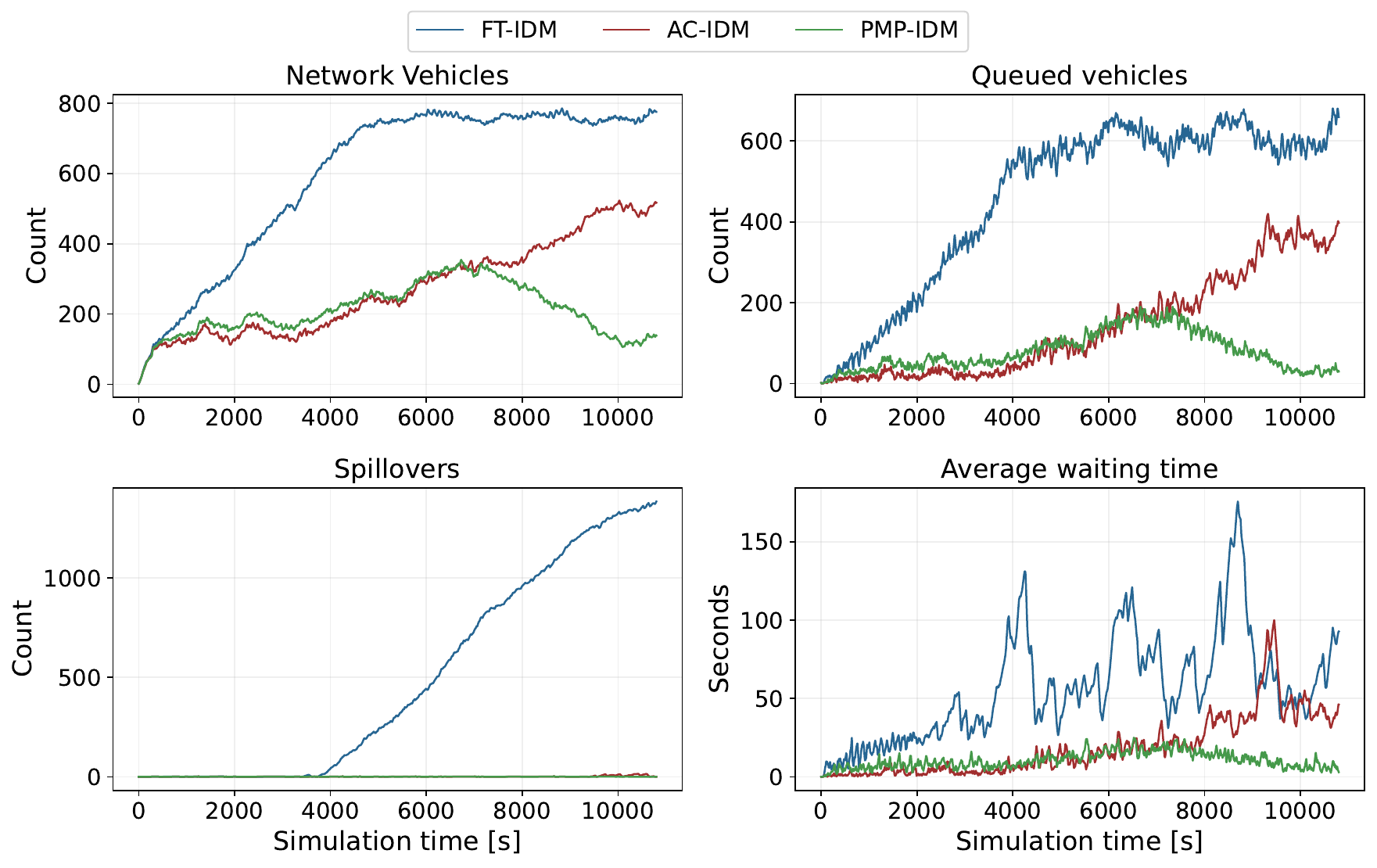}
    \vspace{-0.5cm}
    \caption{Network-level evolution under different signal control methods.}
    \label{fig:ft_2x2}
\end{figure}

In the test stage, the trained SAC agent is deployed to all EBs operating on the eight bus lines in the network. We first examine network-level traffic efficiency. Fig.~\ref{fig:ft_2x2} compares the temporal evolution of (i) the number of vehicles in the network, (ii) queued vehicles, (iii) spillovers, and (iv) average waiting time under different signal control schemes. Under fixed-time control, the network rapidly enters a congested regime: network accumulation and queued vehicles grow monotonically, spillovers emerge after the mid-horizon and continue to increase, and the average waiting time becomes much higher and more volatile, indicating severe gridlock.
In contrast, both actuated control and PMP mitigate congestion relative to fixed-time control, reducing queues and largely suppressing spillovers. Their behaviors, however, differ under sustained high demand. Actuated control performs well in the early stage but gradually deteriorates later, as shown by increasing network accumulation, queue growth, and rising waiting time. By comparison, PMP maintains a more stable evolution: queues and waiting time remain lower and even decline toward the end of the simulation, while spillovers are essentially eliminated. These results indicate that adaptive signal control is critical for this corridor and that pressure-based priority control provides more robust congestion containment than standard actuated control.

\begin{table*}[!t]
\centering
\caption{Network-level summary metrics}
\label{tab:network_summary}

\begin{tabular}{lrrrrrrr}
\toprule
Method
& Car delay [s]
& EB delay [s]
& EB energy {[kWh]}
& EB stops unsched.
& Max vehicles
& Max queued
& Max spillo. \\
\midrule
FT-IDM & 1616.48 & 582.48 & 195.87 & 306 & 787 & 684 & 1386 \\
FT-G2 & 1610.23 & 589.79 & 191.81 & 302 & 782 & 678 & 1367 \\
AC-IDM & 414.26 & 317.96 & 197.80 & 136 & 524 & 421 & 16 \\
AC-G2  & 486.93 & 385.69 & 189.96 & 142 & 602 & 476 & 93 \\
PMP-IDM & 321.03 & 280.72 & 196.51 & 69 & 355 & 195 & 7 \\
PMP-G2 & 347.81 & 323.20 & 192.55 & 68 & 365 & 216 & 7 \\
PriEco-DRL-T & \textbf{292.29} & \textbf{190.13} & 246.57 & 52 & \textbf{316} & \textbf{175} & 7 \\
PriEco-DRL-E & 307.86 & 238.03 & \textbf{185.47} & 49 & 328 & 190 & 7 \\
PriEco-DRL-B & 309.24 & 220.30 & 215.46 & \textbf{48} & 332 & 187 & 7 \\
\bottomrule
\end{tabular}

\end{table*}
Table~\ref{tab:network_summary} summarizes network-level performance in terms of car delay, EB delay, EB energy consumption, and EB unscheduled stops. As expected, the fixed-time baselines (FT-IDM and FT-G2) result in the highest delays for both cars and EBs, as well as severe spillovers, indicating network-wide congestion. Introducing the rule-based G2 strategy under fixed-time signals (FT-G2) yields only marginal improvements in energy consumption and delays, suggesting that when queues and forced stops dominate, vehicle-side speed guidance alone has limited impact without adaptive signal control.

In contrast, adaptive signal control significantly improves traffic efficiency. Both AC and PMP reduce delays and suppress spillovers by responding to real-time traffic conditions. However, the effectiveness of the G2 module depends on the signal controller. Under AC, adding G2 (AC-G2) does not improve delays and increases spillovers, indicating that the speed and dwell advice from G2 is difficult to implement when phase switching is event-driven and unpredictable, particularly under mixed-traffic conditions. Under PMP, PMP-G2 achieves similar spillover suppression as PMP-IDM and yields comparable unscheduled stops, with only marginal changes in energy and delay. This suggests that when a transit-priority adaptive controller already provides favorable progression, adding rule-based guidance offers limited additional benefit.

The three reward-weight settings in PriEco-DRL highlight a clear energy-reliability trade-off. The travel-time-prioritized variant PriEco-DRL-T achieves the lowest EB delay, while the energy-prioritized variant PriEco-DRL-E reduces EB energy consumption the most. The balanced variant PriEco-DRL-B minimizes unscheduled stops, effectively reducing stop-start events under transit-priority signals. Compared to FT-IDM, all PriEco-DRL variants substantially reduce EB delay and avoid network spillovers, offering configurable operating points between energy efficiency and service reliability.

\begin{table*}[t]
\centering
\caption{Ablation study on stage-wise efficiency reward design}
\label{tab:ablation}

\begin{tabular}{lrrrrrrr}
\toprule
Method
& Car delay [s]
& EB delay [s]
& EB energy {[kWh]}
& EB stops unsched.
& Max vehicles
& Max queued
& Max spillo. \\
\midrule

PriEco-DRL-B & \textbf{309.24} & \textbf{220.30} & \textbf{215.46} & \textbf{48} & \textbf{332} & \textbf{187} & 7 \\
PriEco-DRL-B-A1 & 314.79 & 228.31 & 216.93 & 62 & 334 & 187 & 7 \\
PriEco-DRL-B-A2 & 452.38 & 480.15 & 120.56 & 112 & 431 & 267 & 7 \\
\bottomrule
\end{tabular}

\end{table*}

\subsubsection{Ablation study on reward design}
\label{sec:ablation}

To validate the necessity of the proposed stage-wise efficiency reward design, we conduct an ablation study based on the balanced setting PriEco-DRL-B. Table~\ref{tab:ablation} reports network-level metrics for (i) the full model PriEco-DRL-B, (ii) PriEco-DRL-B-A1, which removes the Stage-4 (intersection-crossing) coordination reward while keeping all other terms unchanged, and (iii) PriEco-DRL-B-A2, which removes all stage-dependent efficiency rewards (Stage 1/2/4), such that the policy is driven only by the global energy and safety/comfort terms.

Removing the Stage-4 reward (A1) leads to a clear deterioration in operational reliability: EB unscheduled stops increase from 48 to 62 and EB delay rises moderately, while total EB energy remains nearly unchanged. This confirms that the Stage-4 reward primarily contributes to signal--vehicle coordination around intersections by reducing non-scheduled stop--start events, rather than directly driving energy savings.

In contrast, removing all stage-wise efficiency shaping (A2) causes a substantial degradation in traffic performance, with sharply increased car and EB delays, more unscheduled stops, and larger network accumulation (higher maximum vehicle and queue counts). Although A2 yields lower total energy consumption, this reduction is achieved by an overly conservative low-speed policy that sacrifices service reliability and progression. It therefore represents a degenerate energy-minimizing solution rather than a practical eco-driving strategy.

Overall, the ablation results demonstrate that stage-wise reward shaping is essential for preventing trivial behaviors and for achieving operationally meaningful trade-offs between energy efficiency and travel performance under adaptive, transit-priority signals. The following sections provide more detailed analyses of EB performance at the line and trajectory levels.

\subsubsection{Line-level EB performance}

We next examine line-level EB performance on both the routes encountered during training and the unseen routes used for transfer tests. FT baselines are excluded because gridlock under fixed-time control dominates both travel time and energy use, making the resulting EB metrics less informative for eco-driving comparison. We also omit the G2 variants of AC and PMP here, since their gains under adaptive SPaT are limited and inconsistent at the network level.

\begin{table*}[t]
\centering
\caption{Line-level EB \textbf{energy consumption} comparison across methods. The metrics are averaged over all buses in each line; $\Delta$ is the relative increase/decrease of PriEco-DRL methods compared to the baseline method AC-IDM.}
\label{tab:line_level_energy}
\begin{tabular}{lrrrrrrrrr}
\toprule
\multirow{2}{*}{Line} & \multirow{2}{*}{\parbox[c]{1.4cm}{\centering Travel\\distance [km]}} & \multicolumn{8}{c}{Energy consumption [kWh]} \\
\cmidrule(lr){3-10}
 & & AC-IDM & PMP-IDM & PriEco-DRL-T & $\Delta$ [\%] & PriEco-DRL-E & $\Delta$ [\%] & PriEco-DRL-B & $\Delta$ [\%] \\
\midrule
Line 22-SW & 0.976 & 1.323 & 1.298 & 1.695 & +28.2\% & 1.461 & +10.4\% & 1.594 & +20.5\%\\
Line 22-WN & 0.948 & 1.382 & 1.395 & 1.705 & +23.4\% & 1.309 & -5.3\% & 1.591 & +15.1\%\\
Line 37-W & 2.851 & 4.373 & 4.207 & 5.421 & +24.0\% & 3.734 & -14.6\% & 4.575 & +4.6\%\\
Line 37-E & 2.826 & 4.132 & 4.193 & 5.522 & +33.6\% & 3.674 & -11.1\% & 4.169 & +0.9\%\\
Line 40-SW & 0.805 & 1.116 & 1.184 & 1.300 & +16.4\% & 1.039 & -6.9\% & 1.265 & +13.3\%\\
Line 40-WS & 0.804 & 1.239 & 1.219 & 1.284 & +3.6\% & 1.076 & -13.2\% & 1.210 & -2.3\%\\
Line 65-N & 0.815 & 0.981 & 1.015 & 1.294 & +31.9\% & 1.143 & +16.5\% & 1.293 & +31.7\%\\
Line 65-S & 0.817 & 1.172 & 1.119 & 1.272 & +8.5\% & 1.165 & -0.6\% & 1.284 & +9.5\% \\
\bottomrule
\end{tabular}
\end{table*}
\begin{table*}[t]
\centering
\caption{Line-level EB \textbf{travel time} comparison across methods. The metrics are averaged over all buses in each line; $\Delta$ is the relative increase/decrease of PriEco-DRL methods compared to the baseline method AC-IDM.}
\label{tab:line_level_time}
\begin{tabular}{lrrrrrrrrr}
\toprule
\multirow{2}{*}{Line} & \multirow{2}{*}{\parbox[c]{1.4cm}{\centering Travel\\distance [km]}} & \multicolumn{8}{c}{Travel time [s]} \\
\cmidrule(lr){3-10}
 & & AC-IDM & PMP-IDM & PriEco-DRL-T & $\Delta$ [\%] & PriEco-DRL-E & $\Delta$ [\%] & PriEco-DRL-B & $\Delta$ [\%] \\
\midrule
Line 22-SW & 0.976 & 357 & 345 & 243 & -32.1\% & 266 & -25.6\% & 258 & -27.9\% \\
Line 22-WN & 0.948 & 351 & 331 & 254 & -27.5\% & 278 & -20.7\% & 272 & -22.4\% \\
Line 37-W & 2.851 & 617 & 550 & 456 & -26.2\% & 549 & -11.0\% & 517 & -16.2\% \\
Line 37-E & 2.826 & 725 & 652 & 528 & -27.2\% & 677 & -6.7\% & 614 & -15.3\% \\
Line 40-SW & 0.805 & 311 & 276 & 194 & -37.7\% & 213 & -31.5\% & 214 & -31.3\%\\
Line 40-WS & 0.804 & 285 & 273 & 188 & -33.9\% & 223 & -21.7\% & 202 & -28.9\% \\
Line 65-N & 0.815 & 342 & 283 & 212 & -38.0\% & 234 & -31.4\% & 218 & -36.3\% \\
Line 65-S & 0.817 & 346 & 289 & 209 & -39.6\% & 237 & -31.5\% & 235 & -32.1\% \\
\bottomrule
\end{tabular}
\end{table*}

Tables~\ref{tab:line_level_energy} and \ref{tab:line_level_time} report, respectively, the average total EB energy consumption (kWh) and travel time (s) per line, averaged over all buses in each line, with AC-IDM as the baseline. Overall, PMP-IDM reduces EB travel time on all lines relative to AC-IDM (Table~\ref{tab:line_level_time}), indicating that priority-weighted max-pressure control improves EB progression through intersections. Its effect on energy is more mixed but remains broadly comparable to AC-IDM across lines (Table~\ref{tab:line_level_energy}).
The three PriEco-DRL variants (T/E/B) expose distinct operating points through reward weighting. Table~\ref{tab:line_level_time} shows that all PriEco-DRL variants reduce EB travel time relative to AC-IDM on all eight lines, including the unseen routes, demonstrating transferability beyond the training routes. Among them, the travel-time--oriented variant PriEco-DRL-T achieves the largest reductions (e.g., $-26.2\%$ and $-27.2\%$ on the longer lines Line 37-W and Line 37-E, and up to $-39.6\%$ on Line 65-S), whereas PriEco-DRL-E and PriEco-DRL-B yield slightly smaller but still consistent improvements.

The energy results in Table~\ref{tab:line_level_energy} reveal a more nuanced picture. PriEco-DRL-E achieves the most favorable energy performance overall and reduces energy consumption on several lines (e.g., Line 22-WN, Line 37-W, Line 37-E, Line 40-SW, and Line 40-WS), whereas PriEco-DRL-T tends to increase total energy use on most lines. PriEco-DRL-B generally lies between T and E, although it may still increase energy on some routes. This pattern suggests that, under transit-priority adaptive signals, faster progression may come at the cost of higher energy use, while stronger energy-oriented control encourages smoother and lower-power driving. At the same time, the magnitude of the energy benefit remains route-dependent, reflecting differences in stop spacing, intersection exposure, and mixed-traffic interactions. For example, on the unseen short route Line 65-N, PriEco-DRL-E does not reduce total energy relative to AC-IDM, likely because the baseline energy is already low and the room for additional savings is limited (see its counterpart, the trained Line 65-S).

Taken together, these line-level results show that PriEco-DRL provides configurable operating points for practice: PriEco-DRL-T prioritizes service efficiency, PriEco-DRL-E prioritizes energy savings where feasible, and PriEco-DRL-B offers an intermediate compromise. These observations motivate the trajectory-level analysis in the next section, where we further examine how speed-profile shaping and intersection coordination produce the observed time--energy patterns.

\begin{figure*}[!t]
    \centering

    \subfloat[Line-level EB energy consumption distribution (per bus, normalized by distance)\label{fig:bus-box-energy}]{
        \includegraphics[width=0.80\textwidth]{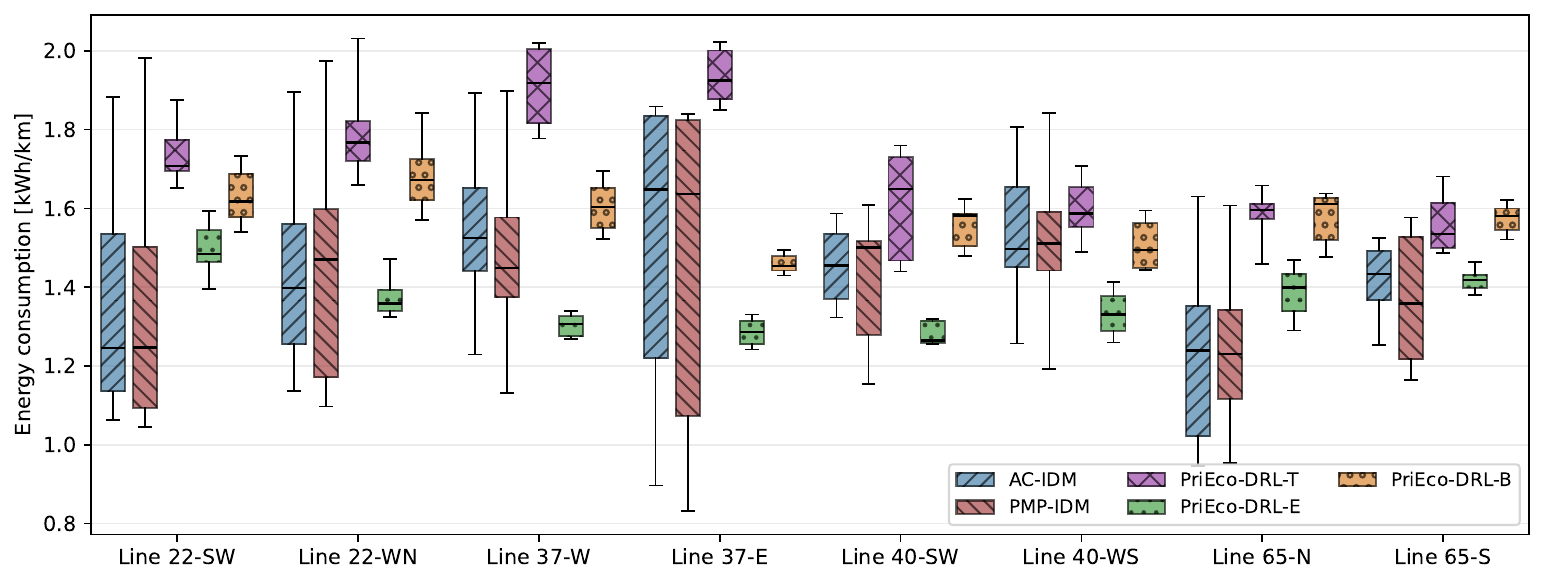}
    }
    

    \subfloat[Line-level EB travel time distribution (per bus, normalized by distance)\label{fig:bus-box-time}]{
        \includegraphics[width=0.80\textwidth]{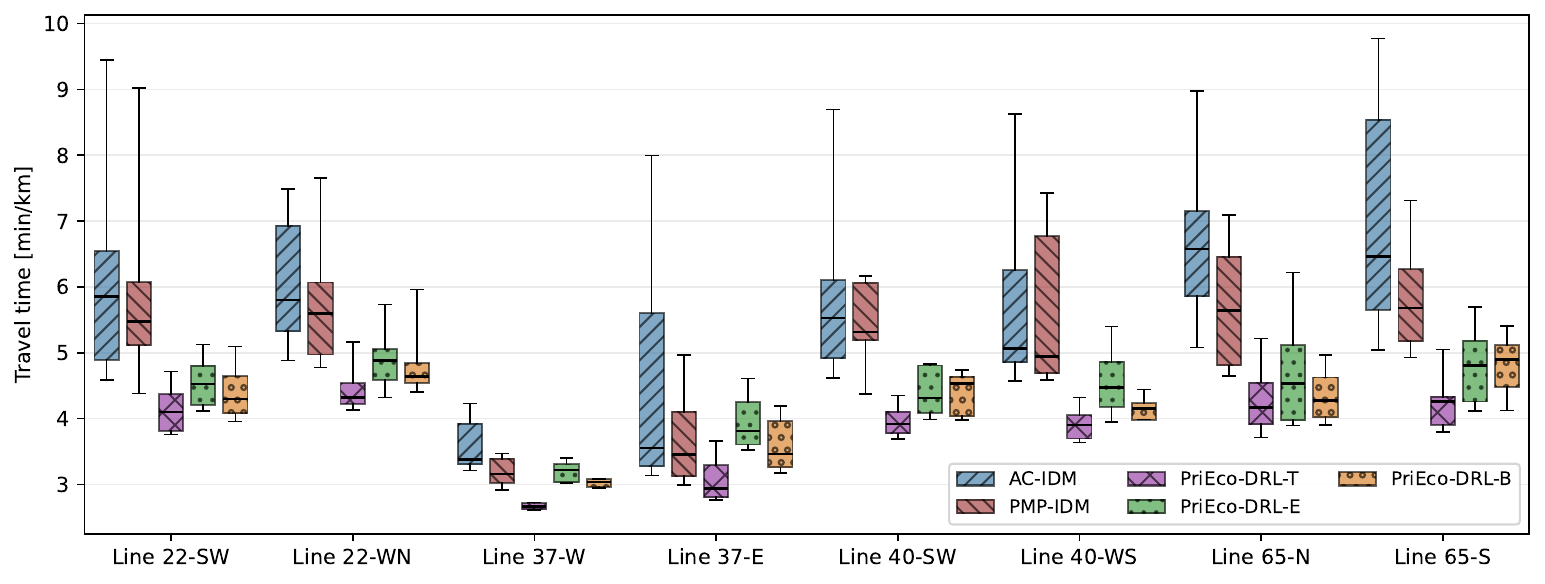}
    }

    \caption{Line-level distributions of individual EB performance across the eight bus lines under different control methods: (a) energy consumption (kWh/km) and (b) travel time (min/km), both normalized by traveled distance.}
    \label{fig:bus-box-plot}
\end{figure*}

Fig.~\ref{fig:bus-box-plot} complements Tables~\ref{tab:line_level_energy}--\ref{tab:line_level_time} by visualizing the distribution of individual-bus outcomes within each line, revealing within-line heterogeneity that is masked by line-averaged totals. To ensure comparability across routes of different lengths, both energy and travel time are normalized by traveled distance (kWh/km and min/km). The travel-time distributions show noticeable spreads for most lines, reflecting varying network conditions, particularly during periods of heavier congestion, when signal delays, residual queues, and mixed-traffic interactions become more significant.

The PriEco-DRL variants distinctly separate the travel-time distributions (panel~b). The time-oriented PriEco-DRL-T consistently shifts the distribution downward, achieving faster progression (lower min/km) than both AC-IDM and PMP-IDM. In contrast, PriEco-DRL-E and PriEco-DRL-B generally show higher min/km, reflecting their stronger focus on energy efficiency or balanced operation. This distribution-level separation visually confirms that reward weighting governs the trade-off between travel efficiency and energy efficiency.

Energy outcomes (panel~a) reveal a more heterogeneous pattern. PriEco-DRL-E tends to lower kWh/km on several lines and often achieves the lowest median energy, indicating that prioritizing energy efficiency drives the policy toward smoother transients and lower-power cruising. Meanwhile, PriEco-DRL-T frequently shows higher kWh/km, consistent with faster operations that require more traction power and less time for speed smoothing. The balanced PriEco-DRL-B generally lies between T and E in both energy and travel time, offering an intermediate behavior.

The PriEco-DRL variants also exhibit tighter distributions on almost all lines, particularly in travel time (panel~b), as indicated by smaller interquartile ranges and shorter whiskers. This suggests that the learned policy improves not only median performance but also \textit{consistency} across individual trips, reducing variability caused by stochastic signal evolution and queue interactions. By regulating approach speeds based on local, pressure-informed signal cues, PriEco-DRL minimizes stop–start perturbations and dampens small timing mismatches, leading to more reliable trip outcomes.

Overall, Fig.~\ref{fig:bus-box-plot} shows that reward weighting controls the operating point in the energy–time space and highlights that PriEco-DRL improves not just average performance but also the \textit{reliability} of EB operations at the individual-bus level.

\subsubsection{Trajectory-level evidence}

Representative trip-level profiles further support the line-level time--energy trade-off observed above. In general, the time-prioritized variant tends to complete the trip earlier at the expense of higher cumulative energy, whereas the energy-prioritized variant favors smoother operation and lower energy use, with the balanced variant lying in between. These trajectory-level patterns are consistent with the aggregate findings. Detailed space--time trajectories and trip-level energy profiles for representative buses are provided in the Supplementary Material.

\section{Conclusion}\label{sec:conclusion}
This paper presents PriEco-DRL, a framework that integrates priority-aware adaptive signal control with DRL-based eco-driving for electric buses (EBs), enabling joint optimization of intersection operations and EB control. The Priority Max-Pressure (Priority-MP) controller dynamically allocates green time, while the DRL policy uses pressure-informed cues for anticipative speed decisions. The CTDE approach with parameter sharing enhances data efficiency and supports deployment across different routes and intersections.

Experiments on a real-world Amsterdam corridor demonstrate that PriEco-DRL significantly reduces EB energy consumption, while maintaining transit priority and avoiding network instability. The energy-prioritized variant achieves the lowest energy use, while the time-prioritized and balanced variants show smaller travel-time increases. The results highlight two key mechanisms: (i) reduced non-scheduled stops due to improved signal–vehicle coordination, and (ii) smoother speed regulation closer to energy-efficient regimes. Ablation analysis confirms that stage-wise reward shaping is crucial for achieving effective energy–reliability trade-offs.

Future work will extend PriEco-DRL to a multi-objective formulation that adjusts the energy–time trade-off online based on changing traffic and EB states, reducing the need for manual reward-weight tuning and offline retraining.

\bibliographystyle{IEEEtran}
\bibliography{IEEEabrv,reference}
\end{document}


\maketitle

\noindent
This supplementary document provides additional experimental details and extended results that support the findings reported in the main manuscript. In particular, it includes additional trajectory-level comparisons, trip-level energy profiles, and representative bus longitudinal profiles that could not be fully presented in the main text due to space limitations.

\section{Additional Trajectory-Level Analyses}
\label{sec:supp_traj}

This section provides additional trajectory-level evidence that complements the line-level results in the main paper. In particular, we present representative space--time trajectories, trip-level energy profiles, and detailed longitudinal speed/acceleration/jerk comparisons for selected EBs from unseen lines. These results further illustrate how PriEco-DRL reduces non-scheduled stop--start events, smooths speed transitions, and selects different operating points along the energy--time trade-off.

\subsection{Representative space--time trajectories}

\begin{figure}[!t]
    \centering

    \begin{subfigure}[t]{1.0\textwidth}
        \centering
        \includegraphics[width=\linewidth, trim={0cm, 0cm, 0cm, 0cm}, clip]{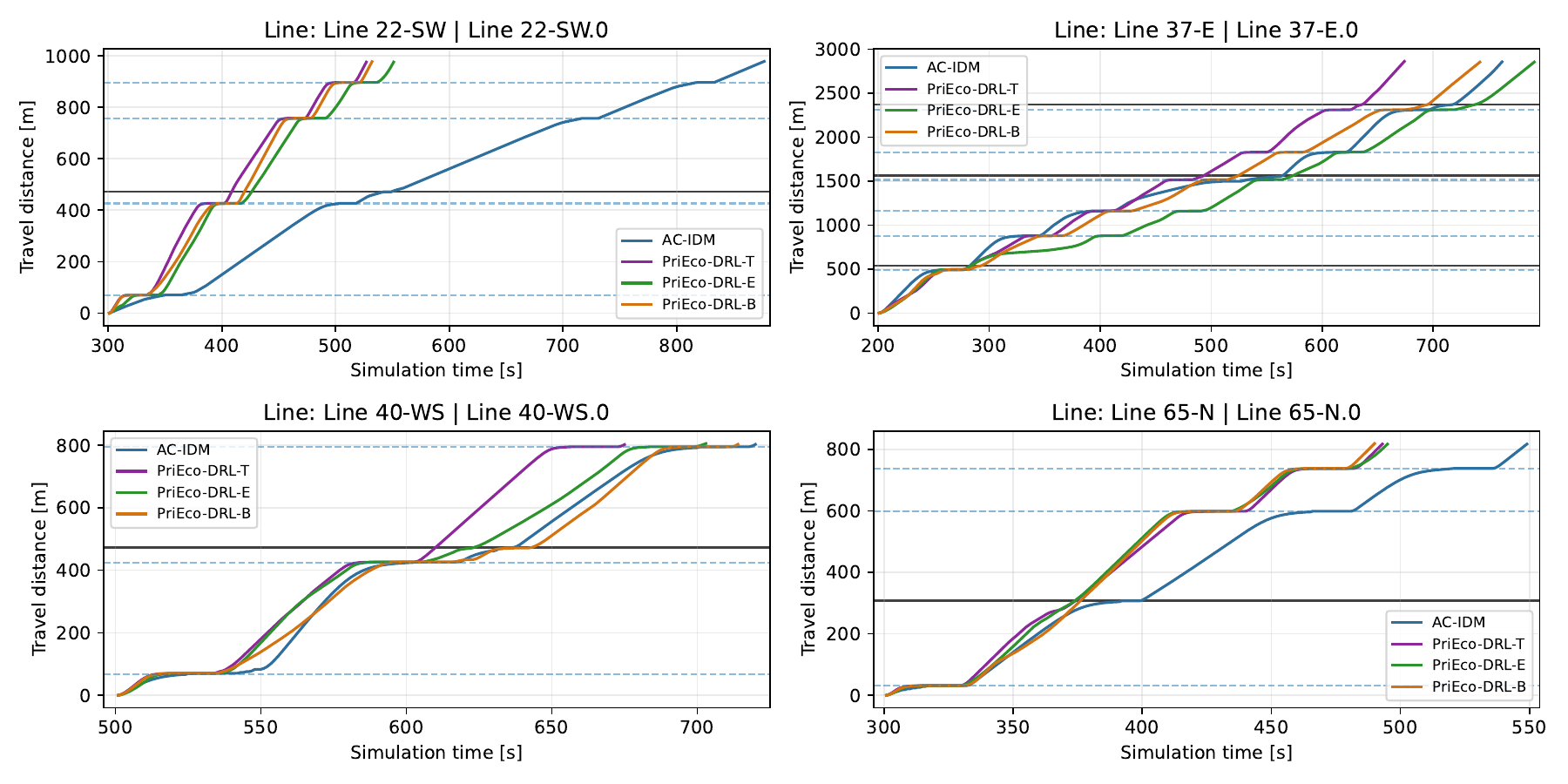}
        \caption{EB trajectories of the unseen four lines at the start of the simulation with light traffic volume}
        \label{fig:supp_bus_traj_start}
    \end{subfigure}

    \vspace{0.6em}

    \begin{subfigure}[t]{1.0\textwidth}
        \centering
        \includegraphics[width=\linewidth]{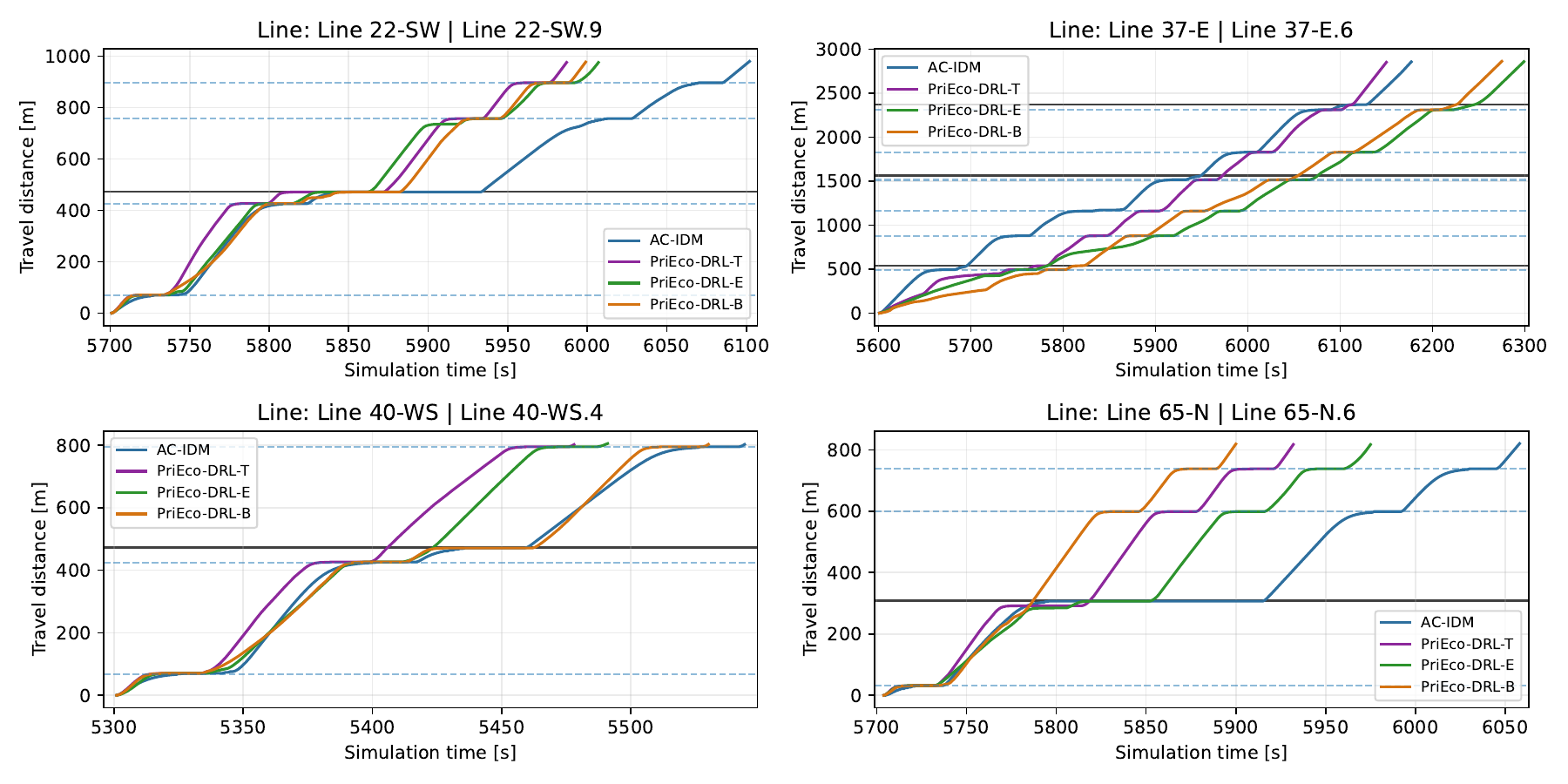}
        \caption{EB trajectories of the unseen four lines in the middle of the simulation with heavy traffic volume}
        \label{fig:supp_bus_traj_mid}
    \end{subfigure}

    \caption{The EB longitudinal space–time trajectories of all the unseen lines under the selected control methods. The solid horizontal line indicates the location of intersection, the dashed blue lines indicate the location of bus stops, and the simulation time indicates the absolute time within the entire simulation.}
    \label{fig:supp_bus_traj}
\end{figure}

Fig.~\ref{fig:supp_bus_traj} presents representative longitudinal space--time trajectories of electric buses from four unseen lines under the selected control methods. Panels (a) and (b) correspond to the beginning and the middle of the 3-hour simulation, respectively. Because adaptive signal controllers generate method-specific phase sequences, the exact switching instants are not directly comparable across methods; therefore, only intersection locations (solid horizontal lines) and bus stops (dashed lines) are marked.

During the low-demand period (Fig.~\ref{fig:supp_bus_traj_start}), AC-IDM exhibits clear progression slowdowns around intersections, with longer flat segments near the intersection markers indicating non-scheduled waiting. In contrast, all PriEco-DRL variants achieve faster progression and shorter non-scheduled holds, which is consistent with the line-level reduction in travel time reported in the main text. The difference is especially visible on longer trajectories such as Line 37-E.0, where PriEco-DRL maintains more continuous advancement over successive intersections. Among the variants, PriEco-DRL-T generally exhibits the steepest effective slope, reflecting a stronger emphasis on travel-time efficiency, whereas PriEco-DRL-E adopts a more conservative progression pattern. PriEco-DRL-B typically lies between the two.

In the congested period (Fig.~\ref{fig:supp_bus_traj_mid}), all methods experience slower progression due to longer queues and intermittent spillbacks, and some buses (e.g., Line 22-SW.9) undergo extended waiting at intersections. Under such conditions, uninterrupted green access becomes difficult for all methods, so the main advantage of PriEco-DRL comes from bus-side pacing rather than guaranteed non-stop passage. Even in these cases, PriEco-DRL reduces large stop--start disruptions and yields a more regular advancement pattern. For example, for Line 37-E.6 and bus\_40Mui.4, the PriEco-DRL trajectories progress more steadily across successive intersections and show shorter non-scheduled holds than AC-IDM. Overall, these trajectory-level observations support the two main mechanisms identified in the paper: (i) reducing intersection-induced stop--start events through anticipative approach-speed regulation, and (ii) shaping smoother speed profiles that avoid unnecessary high-power transients, with the specific pacing determined by the reward-weight setting.

\subsection{Additional trip-level energy profiles}

\begin{figure*}[!t]
    \centering
    \includegraphics[width=\textwidth]{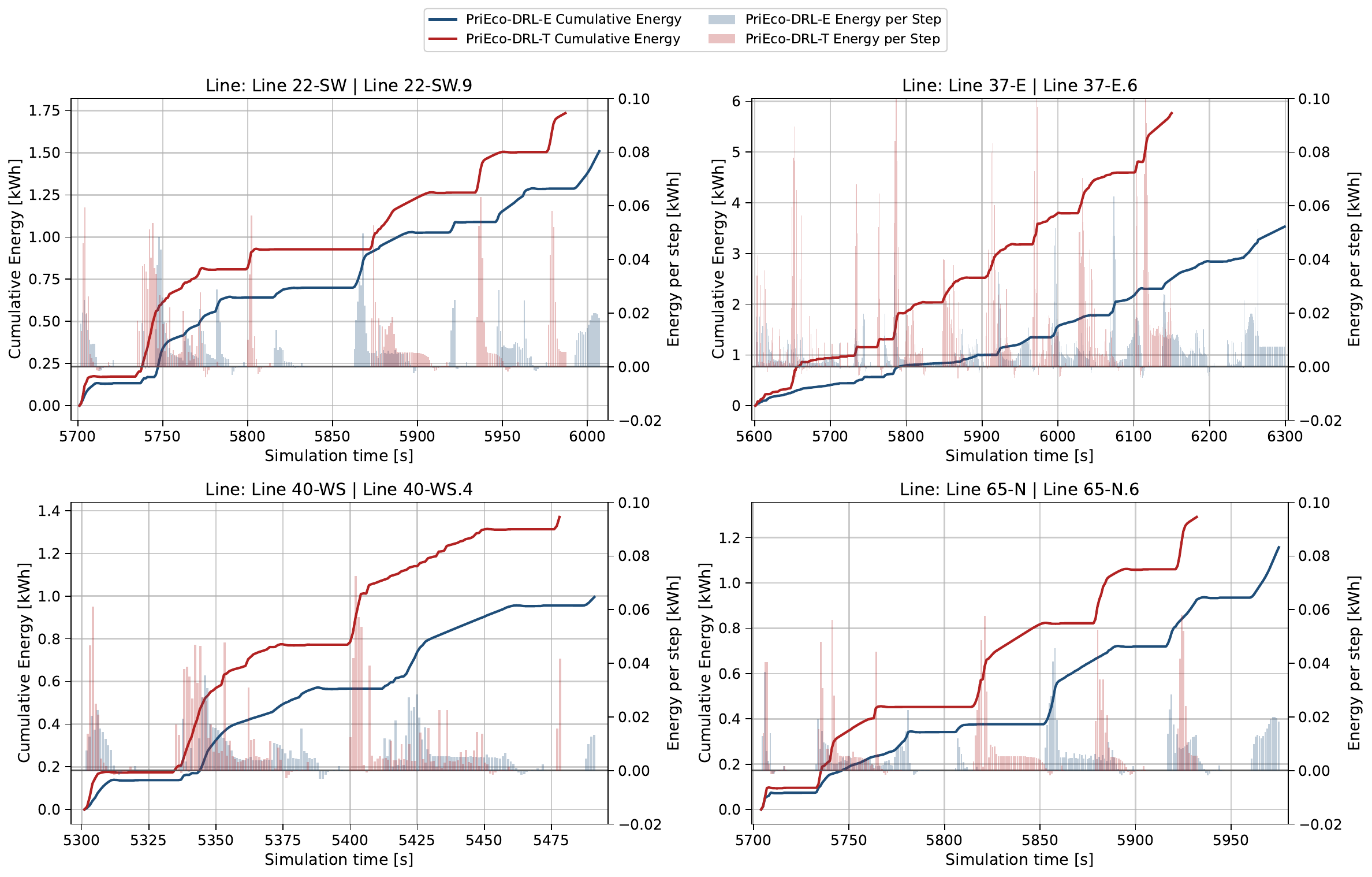}
    \caption{Representative cumulative electrical energy trajectories and step-wise net energy increments for unseen lines during the congested period.}
    \label{fig:supp_energy}
\end{figure*}

Fig.~\ref{fig:supp_energy} further illustrates the trip-level energy-consumption dynamics by plotting (i) the cumulative electrical energy consumption (solid curves, left axis) and (ii) the step-wise net energy increment (shaded bars, right axis) for four representative EBs from unseen lines during the congested period. The shaded series reveals how energy accumulates under mixed traffic and adaptive SPaT: large positive spikes typically correspond to high-power traction events such as re-acceleration after queuing or stop--go recovery, whereas smaller increments indicate smoother progression with fewer abrupt power demands. Negative increments, when present, correspond to regenerative braking, i.e., partial recovery of kinetic or potential energy during deceleration.

Comparing the two reward-weight variants, PriEco-DRL-E consistently yields lower cumulative energy trajectories than PriEco-DRL-T across all four trips, and the gap grows as the trip proceeds. This difference can be traced to the step-wise profiles: PriEco-DRL-T exhibits more frequent and larger energy spikes (e.g., for Line 37-E.6 and Line 22-SW.9), indicating repeated high-power transients associated with aggressive progression and stop--go recovery in congestion. In contrast, PriEco-DRL-E suppresses both the frequency and magnitude of these spikes, suggesting that stronger energy-oriented reward weighting encourages more moderate approach and discharge speeds, fewer unnecessary stop--start maneuvers, and a lower-variance speed profile. These effects accumulate over time and explain the lower total energy use of PriEco-DRL-E observed in the line-level results.

\subsection{Representative EB longitudinal profiles}

\begin{figure*}[!t]
    \centering
    \includegraphics[width=\textwidth]{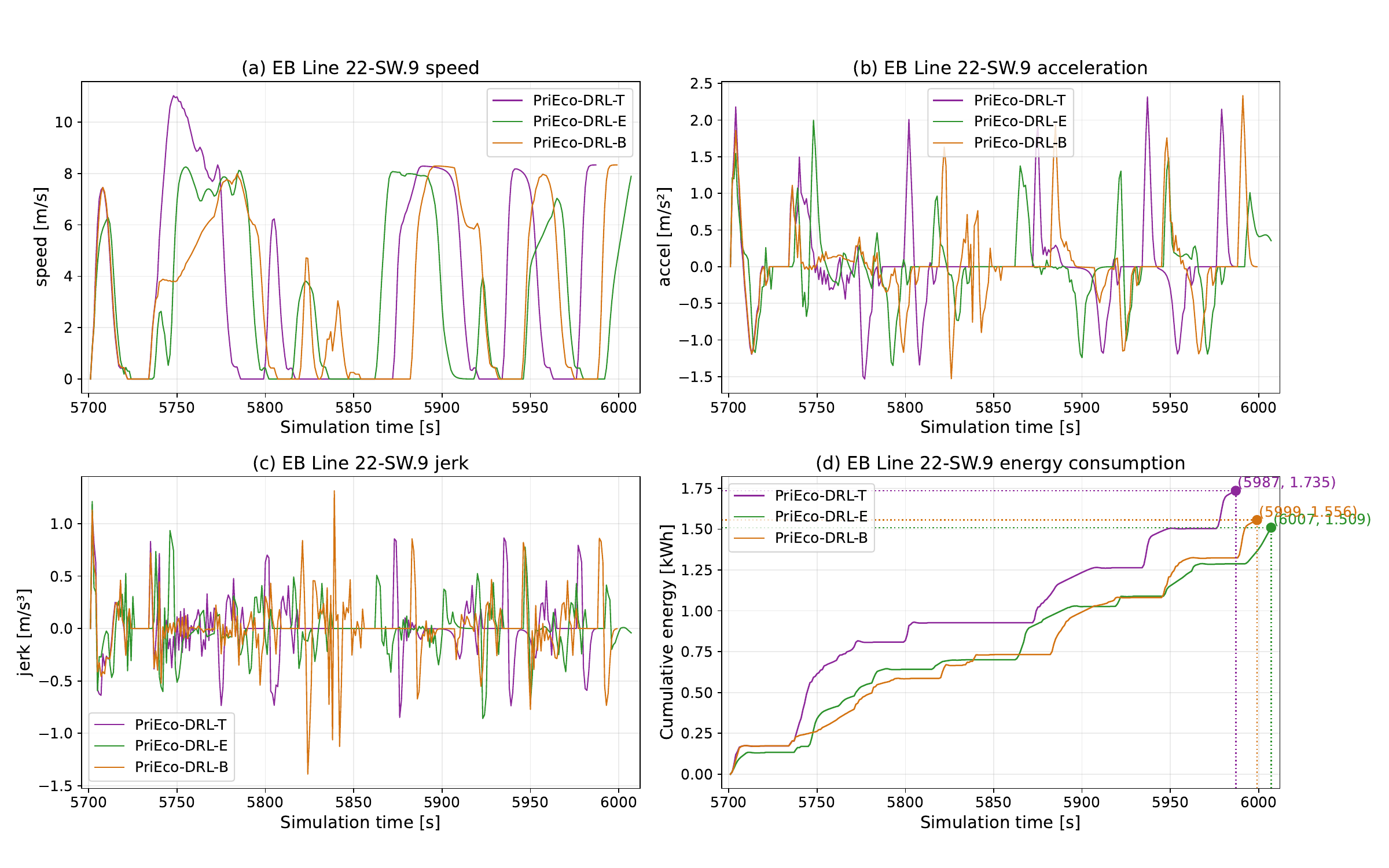}
    \caption{Representative longitudinal profiles of Line 22-SW.9 under PriEco-DRL-T, PriEco-DRL-E, and PriEco-DRL-B during the congested period: (a) speed, (b) acceleration, (c) jerk, and (d) cumulative electrical energy consumption.}
    \label{fig:supp_bus_profile}
\end{figure*}

To further explain the time--energy trade-off among the PriEco-DRL variants, Fig.~\ref{fig:supp_bus_profile} compares the longitudinal profiles of a representative EB trip (Line 22-SW.9) during the congested period. Panels (a)--(c) show the speed, acceleration, and jerk trajectories, and panel (d) reports cumulative energy consumption.

The time-prioritized variant PriEco-DRL-T pursues more aggressive progression, reaching higher cruising speeds and exhibiting more pronounced acceleration events. In contrast, the energy-prioritized variant PriEco-DRL-E adopts a smoother and more conservative pacing strategy, with fewer high-speed excursions and milder acceleration/deceleration swings. PriEco-DRL-B provides an intermediate behavior, combining relatively steady progression with moderated actuation, although localized transients can still appear depending on the surrounding traffic and signal interactions.

These behavioral differences are directly reflected in the cumulative energy trajectories. PriEco-DRL-T completes the trip earlier but consumes the most energy, PriEco-DRL-E achieves the lowest energy at the cost of a longer completion time, and PriEco-DRL-B lies in between. The three trip-end points therefore form a clear Pareto-like ordering for this representative trip, which is consistent with the system-level energy--time trade-off discussed in the main paper. Overall, the example illustrates two key contributors to energy reduction in PriEco-DRL: (i) suppressing energy-intensive stop--start recovery through anticipative speed regulation, and (ii) avoiding unnecessary high-speed cruising and large actuation transients so that the bus operates closer to an energy-efficient regime.

